\documentclass{amsart}
\usepackage{graphicx}
\usepackage{epsfig,latexsym,amsfonts,amssymb,amsmath,amscd,graphics,epic}
\usepackage{amsfonts,amssymb,amsmath,amscd,amsthm}
\usepackage[mathscr]{eucal}
\usepackage{mathrsfs}

\usepackage{mathbbol}

\usepackage{oldgerm,units}
\usepackage{wrapfig,epsfig}

\usepackage{ifthen}
\usepackage{mathbbol}

\usepackage{amsthm}
\usepackage[mathscr]{eucal}
\usepackage{mathrsfs}
\usepackage{mathbbol}
\usepackage{oldgerm,units}
\usepackage{wrapfig}

\newtheorem{theorem}{Theorem}[section]
\newtheorem{proposition}[theorem]{Proposition}
\newtheorem{definition}[theorem]{Definition}

\newtheorem{lemma}[theorem]{Lemma}

\newtheorem{corollary}[theorem]{Corollary}

\newtheorem{example}[theorem]{Example}
\newtheorem{remark}[theorem]{Remark}

\newtheorem{note}[theorem]{Note}

\newtheorem*{note*}{Note}

\newcommand{\Ref}[1]{(\ref{#1})}

\newcommand{\Real}{\mathbb R}

\newcommand{\Int}{\mathbb Z}
\newcommand{\Net}{\mathbb N}

\newcommand{\Trop}{\mathbb T}

\newcommand{\eInt}{\bar{\Int}}

\newcommand{\trop}[1]{\mathcal{#1}}

\newcommand{\tB}{\trop{B}}

\newcommand{\tD}{\trop{D}}

\newcommand{\tF}{\trop{F}}

\newcommand{\tM}{\trop{M}}
\newcommand{\tN}{\trop{N}}

\newcommand{\tS}{\trop{S}}

\newcommand{\tU}{\trop{U}}

\newcommand{\tW}{\trop{W}}

\newcommand{\tpF}{f}

\newcommand{\tfF}{\tilde{f}}

\newcommand{\To}{\longrightarrow }

\newcommand{\mTo}{\longmapsto}

\newcommand{\tUniS}{-\infty}

\newcommand{\al}{\alpha}
\newcommand{\bt}{\beta}

\newcommand{\sig}{\sigma}

\newcommand{\Gm}{\Gamma}

\newcommand{\lm}{\lambda}
\newcommand{\Lm}{\Lambda}

\newcommand{\TrS}{\oplus}
\newcommand{\TrP}{\odot}

\newcommand{\OP}{\left(}
\newcommand{\CP}{\right)}

\def \aaa{a}

    \ifx\proof\undefined
    \newenvironment{proof}{
    \smallskip
    \noindent\emph{Proof.}}{\hfill\(\Box\)
    \bigskip
    } \fi

\newcommand{\vMat}[4]{\OP \begin{array}{cc}
  #1 & #2 \\
  #3 & #4
\end{array}\CP}

\newcommand{\vvMat}[9]{\small{\OP \begin{array}{ccc}
  #1 & #2 & #3\\
  #4 & #5 & #6\\
  #7 & #8 & #9\\
\end{array}\CP}}

\newcommand{\bfem}[1]{\textbf{\emph{#1}}}

\newcommand{\ifdef}[3]{\ifthenelse{\equal{#1}{true}}{#2}{#3}}

\def\pSkip{\vskip 1.5mm \noindent}

\newcommand{\etype}[1]{\renewcommand{\labelenumi}{(#1{enumi})}}

\def\eroman{\etype{\roman}}

\def\mInf{\tUniS}

\def\mT{M_2(\Trop)}
\def\uT{U_2(\Trop)}
\def\mpT{M_2(\Trop[\lm_1, \dots, \lm_6])}
\def\upT{U_2(\Trop[\lm_1, \dots, \lm_6])}

\def\tf{{\tilde f}}
\def\bool{\mathbb B}

\def\lT{L_2(\Trop)}

\def\rep{R}

\def\taA{A}
\def\taB{B}
\def\taI{I}
\def\taZ{Z}
\def\itaA{\taA^{\nabla}}

\def\eqR{\overset{_e}{\sim}}
\def\neqR{\overset{_e}{\nsim}}

\def\l{\operatorname{\ell}}
\def\r{\operatorname{r}}
\newcommand{\mtrace}[1]{\operatorname{Tr_{\TrP}}(#1)}

\newcommand{\mat}[1]{{#1}^{\operatorname{mtx}}}
\newcommand{\pol}[1]{{#1}_{\operatorname{ply}}}

\newcommand{\per}[1]{|{#1}|}
\newcommand{\Per}[1]{\left|{#1}\right|}
\newcommand{\rank}[1]{\operatorname{rank}({#1})}
\newcommand{\adj}[1]{\operatorname{Adj}({#1})}
\newcommand{\trn}[1]{{#1}^{\operatorname{t}}}

\def\sHom{\varphi}
\def\pmHom{\mu}

\def\seq{\wp}

\def\ef{f^e}

\def\eg{g^e}

\def\bfA{\textbf{A}}

\def\bfa{\textbf{a}}

\def\bfi{ \textbf{i}}
\def\bfj{\textbf{j}}

\def\pA{\pol{A}}
\def\pB{\pol{B}}
\def\pF{\pol{F}}
\def\pG{\pol{G}}

\begin{document}


\title[Semigroup Identities in the Monoid of Two-by-Two Tropical Matrices]
{Semigroup Identities in the Monoid \\ \vskip 1mm  of Two-by-Two
Tropical Matrices}


\author{Zur Izhakian }
\thanks{The first author has been supported by the
Chateaubriand scientific post-doctorate fellowships, Ministry of
Science, French Government, 2007-2008; and was partially supported
by a grant from the European Science Foundation (ESF), Automata:
from Mathematics to Applications, No. 1609, 2007}

\address{Department of Mathematics, Bar-Ilan University, Ramat-Gan 52900,
Israel} \address{ \vskip -6mm CNRS et Universit´e Denis Diderot
(Paris 7), 175, rue du Chevaleret 75013 Paris, France}
\email{zzur@math.biu.ac.il, zzur@post.tau.ac.il}
\author{Stuart W. Margolis}
\address{Department of Mathematics, Bar-Ilan University, Ramat-Gan 52900,
Israel}\email{margolis@math.biu.ac.il}


\date{February 2008}
\date{\today}

\keywords{Tropical (Max-Plus) Matrix Algebra, Semigroup Identity,
Idempotent Semiring, Monoid Representation}


\begin{abstract} We show that the monoid  $\mT$ of $2 \times 2$ tropical matrices
is a regular semigroup satisfying the semigroup identity
$$ A^2 B^4 A^2  \ A^2 B^2 \  A^2 B^4 A^2   =   A^2 B^4 A^2 \ B^2
A^2 \
 A^2 B^4 A^2 \ .$$
Studying reduced identities for subsemigroups of $\mT$, and
introducing a faithful semigroup representation for the bicyclic
monoid by $2 \times 2$ tropical matrices, we reprove Adjan's
identity for the bicyclic monoid in a much simpler way.
\end{abstract}

\maketitle




\section*{Introduction}

Varieties of semigroups have been intensively studied for many
years. It is known that the group of all invertible $2 \times 2$
matrices over a field of characteristic 0 contains a copy of the
free group and thus does not satisfy any group or semigroup
identities. Thus the monoid of all $2 \times 2$ over this field
generates the variety of all monoids and semigroups and the $2
\times 2$ general linear group generates the variety of all
groups.

In the last years, tropical mathematics, that is mathematics based
upon the tropical semiring, has been intensively studied. In
particular, the monoid and semiring of $n \times n$ matrices
plays, as one would expect, an important role both algebraically
and in applications to combinatorics and geometry. In contrast to
the case of matrices over a field, we identify a non-trivial
semigroup identity satisfied by the monoid $\mT$ of all $2 \times
2$ tropical matrices and for some of its submonoids. We also note
that the group of units of this monoid is virtually Abelian and
thus both the monoid of all $2 \times 2$ tropical matrices and its
group of units generate proper varieties of monoids and groups
respectively.

Tropical mathematics  has been developed mostly over the tropical
semiring $\Trop\ = \Real\cup\{-\infty\}$ with the operations of
maximum and summation,
$$a\TrS b=\max\{a,b\},\quad a\TrP b=a+b,$$
as addition and multiplication respectively
\cite{IMS,pin98,RST}. It is natural for developing the connections
between classes of semigroups and their matrix representations by
considering matrices over $\Trop$  as the target for representing
semigroups.

One of the fundamental properties of a semigroup is being regular.
We prove that $\mT$ is a regular monoid in Von-Neumann's sense,
and indicate a naturally occurring generalized inverse for each
matrix in $\mT$. We present a few semigroup identities for
submonoids of $\mT$ and particularly for $\mT$ itself: \pSkip
 \textbf{Theorem \ref{thm:idTr}}: \emph{The submonoid $\uT \subset \mT$ of
upper triangular tropical matrices satisfies the semigroup
identity
$$ AB^2A \ AB  \ AB^2A = AB^2A \  B A \ AB^2A \ ;$$}

We note that this is precisely Adjan's identity, the identity of
smallest length satisfied by the bicyclic monoid \cite{Adjan}. In
fact, we prove that again in contrast with the case of matrices
over a field, the bicyclic monoid has a faithful representation in
$\uT$. This opens up the possibility of using representation
theory over the tropical semiring to study the bicyclic monoid. In
particular, we use this faithful tropical linear representation
for the bicyclic monoid in $\uT$ to reprove Adjan's identity in a
much shorter and friendlier way.

Using the above, we prove the main result of the paper.

\pSkip

\noindent \textbf{Theorem \ref{thm:globalId}:} \emph{The monoid
$\mT$ admits the semigroup identity
$$A^2B^4A^2 \ A^2B^2  \ A^2B^4A^2 =A^2B^4A^2 \  B^2 A^2 \ A^2B^4A^2
\ .  $$}

In the past years, most of the theory of matrix semigroups has
been developed for matrices built over fields and rings. In this
paper we appeal to matrices built over semirings which we believe
are the ``current'' structure to establish representations for
classes of semigroups. The bicyclic monoid is one main supporting
example.

\section{Tropical Semirings}
We open by reviewing some basic notions of tropical algebra and
geometry, including the corresponding categorical framework, and
introduce new tropical notions which will be used later in our
exposition.


\subsection{Tropical polynomials}

 Elements of
$\Trop[\lm_1,\dots,\lm_m]$ are called tropical polynomials in $m$
variables over $\Trop$, and are of the form
\begin{equation}\label{eq:polyToFunc1} \tpF = \bigoplus_{\bfi \in
\Omega} \al_\bfi  \lm_1^{i_1}  \cdots  \lm_m^{i_m}  \ \in \
\Trop[\lm_1,\dots,\lm_m]\backslash\{-\infty\},
\end{equation}
where $\Omega\subset\Int^{(m)}$ is a finite nonempty set of points
$\bfi =(i_1,\dots,i_m)$ with nonnegative coordinates, $\al_\bfi
\in\Real$ for all $\bfi \in\Omega$, and  $\al^i$ means $\al
\TrP\cdots\TrP \al$ with $\al$ repeated $i$ times. Any tropical
polynomial determines a piecewise linear convex function
$\tfF:\Real^{(m)} \to \Real$, defined by:
\begin{equation}\label{eq:polyToFunc2} \tfF(\bfa)=\max_{\bfi \in\Omega}\{ \langle
\bfa ,\bfi \rangle+ \al_\bfi \} ,\quad  \bfa \in\Real^{(m)}\ ,
\end{equation}
where $\langle \, \cdot , \cdot \, \rangle$ stands for the
standard scalar product. The map $\tpF\mapsto\tfF$ is not
injective and one can reduce the polynomial semiring so as to have
only those elements needed to describe functions.

Given a point $\bfa = (a_1, \dots, a_m) \in \Trop^{(m)}$, there is
a tropical semiring homomorphism
\begin{equation}\label{eq:subsHom} \sHom_\bfa:
 \Trop[\lm_1,\dots, \lm_m] \ \To \  \Trop
\end{equation}
given by sending
 $$\sHom_\bfa:\bigoplus _\bfi \al_{\bfi} {\lm _1}^{i_1}\cdots  {\lm _m}^{i_m}
 \  \longmapsto \
 \bigoplus _\bfi \al_{\bfi} {a _1}^{i_1}\cdots  {a _m}^{i_m}$$
which we call the \textbf{substitution homomorphism} (relative to
$\bfa$). We write $f(\bfa)$ for the image of $f$ under
$\varphi_\bfa$ and identify each  $\al_\bfi \in \Trop$ with the
monomial $\al_\bfi \lm_1^0 \cdots \lm_m^0$ to have the embedding
$\Trop \hookrightarrow \Trop[\lm_1,\dots, \lm_m]$.

\begin{definition}\label{def:essentialPart} Suppose
$f = \bigoplus \al _\bfi   \lm_1^{i_1}  \cdots  \lm_m^{i_m},$ $h =
\al_\bfj   \lm_1^{j_1}  \cdots  \lm_m^{j_m}$ is a monomial of $f$,
and write $f_h = \bigoplus _{\bfi \ne \bfj}\al _\bfi \lm_1^{i_1}
\cdots  \lm_m^{i_m}.$ We say that the monomial $h$ for $f$  is
\textbf{inessential} if $f(\bfa) = f_h(\bfa)$ for each $\bfa \in
\Trop^{(m)}$; otherwise $h$ is said to be \textbf{essential}. The
\textbf{essential part} $\ef$ of a polynomial
 $f = \bigoplus \al _\bfi   \lm_1^{i_1}  \cdots  \lm_m^{i_m}$ is the sum of those monomials
$\al_\bfj   \lm_1^{j_1} \cdots  \lm_m^{j_m}$ that are essential.
When $f = \ef$, $f$ is said to be an \textbf{essential
polynomial}.
\end{definition}
\noindent (Note that,  any  monomial $\neq \mInf$ by itself,
considered as a polynomial,
 is essential.)

A monomial $h$ is essential in a polynomial $f$ if $h(\bfa)>
f_h(\bfa)$ for some $\bfa$ and thus for all $\bfa'$ in an open set
$W_\bfa$ of the standard topology of $\mathbb{R}^{m}$containing
$\bfa$. Any monomial $h$ of $\ef$ is essential in $\ef$. Indeed,
by definition, $f_h(\bfa) \TrS h(\bfa)> f_h(\bfa)$ for some
$\bfa\in \Trop^{(m)},$ implying $ h(\bfa)
> f_h(\bfa)$.

Using this definition we say that two polynomials $f$ and $g$ are
\textbf{essentially equivalent}, written $f \eqR g$, if $\ef =
\eg$. We shall  show that the essential part of a polynomial $f$
defines the same function as $f$, that is, $ f(\bfa) = g(\bfa)$
for each $\bfa \in \Trop^{(m)}$ and is the unique essential
polynomial with this property. Thus $\ef$ is a canonical
representative for the congruence class of the morphism that sends
each polynomial to the function it defines.

\begin{remark}\label{arch1}
For any nonconstant monomials $g_1,g_2,h_1, \dots, h_k$ and
$\bfa\in \Trop^{(m)}$ with
$$g_2(\bfa) = g_1(\bfa )
> h_i(\bfa), \quad 1 \le i \le k,$$ and any open set
$W_\bfa$ of $\Trop^{(m)},$ containing $\bfa$,  there exists
$\bfa'\in W_\bfa$ with
$$g_2(\bold {a'}) > g_1(\bold {a'})
> h_i(\bold {a'}), \quad 1 \le i \le k.$$
\end{remark}

\begin{lemma} For any monomials $g_1,\dots, g_\ell$,
$h_1, \dots, h_k$ and $\bfa \in \Trop^{(m)}$ with
$$g_1(\bfa ) = g_2(\bfa) = \dots =g_\ell(\bfa)
> h_i(\bfa ), \quad 1 \le i \le k,$$ there exists  $\bfa' \in \Trop^{(m)}$
and $1<j \le \ell$ such that $$g_j(\bold {a'})
> g_i(\bold {a'}) \quad \forall i \ne j; \qquad g_j(\bold {a'})
> h_i(\bold {a'}) , \quad 1 \le i \le k.$$
\end{lemma}

\begin{proof} Induction on $\ell$. By Remark \ref{arch1}, we have
$\bfa'\in \Trop^{(m)}$ such that $$g_2(\bfa')> g_1(\bfa ')
> h_i(\bfa'), \quad 1 \le i \le k.$$ Take $j$ such that $g_j(\bfa')$ is maximal,
and expand the $h_i$ to include all $g_i$ such that $g_j(\bfa')>
g_i(\bfa').$ Then we have the same hypothesis as before, but with
smaller $\ell$.
\end{proof}

\begin{proposition}\label{prop:essential}  $\ef$ defines the same function as $f$
for any $f\in \Trop[ \lm_1,\dots,\lm_m] $.
\end{proposition}
\begin{proof} Given any $\bfa  \in \Trop^{(m)}$, there is a monomial
$g_1$ such that $f(\bfa ) = g_1(\bfa ).$ We need to show that
$\ef(\bfa ) = g_1(\bfa ) .$ Suppose  $g_1(\bfa ) = g_2(\bfa ) =
\cdots = g_\ell(\bfa )
> h(\bfa ) $ for some other monomial(s) $g_2, \dots, g_\ell$ of
$f$ which are inessential in $f$. But then, by the lemma, we may
find $\bfa'$ such that $g_j(\bfa')$ takes on the single largest
value of the monomials of $f$, for some $2 \le j \le \ell,$
contrary to $g_j$ being inessential in $f$.
\end{proof}

\begin{corollary} $\ef$ is well defined for any polynomial $f$.
\end{corollary}

\begin{proof} Assume $\ef_1 \neq \ef_2$ are two different
essential parts of $f$. Then by Proposition \ref{prop:essential}
each of them defines the same polynomial function as $f$. Since
$\ef_1 \neq \ef_2$, there exists $\bfa \in \Trop^{(m)}$ such that
$\ef_1(\bfa) = g_1(\bfa) \neq h_1(\bfa) =\ef_2(\bfa)$ for
monomials $g_1$ and $h_1$ of $\ef_1$ and $\ef_2$, respectively.
Suppose $g_1(\bfa) > h_1(\bfa)$, then by Lemma \ref{arch1}, there
is an open set $W \subset \Trop^{(m)}$ containing $\bfa$ such that
$\ef_1(\bfa') =g_1(\bfa')
> h_1(\bfa')= \ef_2(\bfa')$ for each $\bfa' \in W$ that is $\ef_1$ defines a different function from $\ef_2$
-- a contradiction.
\end{proof}

Therefore, the polynomial semiring $\Trop[\lm_1,\dots,\lm_m]$ can
be viewed as the collection of  essential polynomials, viewed as a
semiring where we perform the usual operations in $\Trop[
\lm_1,\dots,\lm_m]$ and then take the essential part. Clearly,
$\eqR$ is an equivalence relation and, for short, we call it
$e$-equivalence. It is easy to see that $\eqR$ is a semiring
congruence, and thus we have the quotient semiring $
\Trop[\lm_1,\dots,\lm_m] / {\eqR} $. Together with Equation
\Ref{eq:polyToFunc2}, this gives the semiring homomorphism
$$ \phi : \Trop[\lm_1,\dots,\lm_m] / {\eqR} \ \To  Poly(\Trop^{m}),$$
where $Poly(\Trop^{m})$ is the semiring of $\Trop$-valued
polynomial functions.


%
%

\begin{lemma} The homomorphism $\phi$ is an isomorphism.
That is, every polynomial function has a unique essential
representative.
\end{lemma}
\begin{proof} Recall that a polynomial function $\tf : \Real^{(m)} \to \Real$ of the form
\eqref{eq:polyToFunc2} determines a piecewise linear convex
function whose graph $\Gm_\tf \subset \Real^{(m+1)}$ consists of a
finite number of facets (i.e. faces of co-dimension $1$). Each
facet $\tF$ of $\Gm_\tf$ is described uniquely by a linear
equation (in the classical sense), that is a tropical monomial.
Accordingly, given a function which is determined by an essential
polynomial $\ef$ with $\ell$ monomials, its graph will have
exactly $\ell$ facets, each of them is the image of a subset
$W_{i} \subset \Real^{(m)}, 1 \leq i \leq l$ on which the
evaluation of $\ef$ is attained by a single monomial which takes
maximal values over all the other monomials of $\ef$ on $W$.

Suppose that there are two essential polynomials $\ef_1$ and
$\ef_2$ corresponding to the same function $\tf$ whose graph
$\Gm_\tf$ has $\ell$ facets. Then $\ef_1$ and $\ef_2$ have exactly
$\ell$ monomials, each monomial corresponds to a facet $\tF$ of
$\Gm_\tf$. Take a facet $\tF$. Since it has a unique description
each pair of monomials of $\ef_1$ and $\ef_2$ corresponding to
$\tF$ must be identical. This is true for each facet, and
therefore  $\ef_1= \ef_2$.
\end{proof}

We say that a polynomial $f$ is a \textbf{flat polynomial} if all
of its  coefficients are equal to some fixed value $\al \in
\Trop^\times$, where tropically $\Trop^\times$ stands for  $\Trop
\setminus \{ \tUniS \}$.

\begin{lemma}\label{lem:essential} Assume $f = f_1 \TrS  f_2 \TrS f_3$
is a flat polynomial in $\Trop[\lm_1,\dots,\lm_m]$, where \pSkip
$$ f_1 = \lm_{1}^{j_1+k}\lm_{2}^{j_2-k} \lm_{3}^{j_3}\cdots
\lm_{m}^{j_{m}}, \qquad f_2 = \lm_{1}^{j_1}\lm_{2}^{j_2}
\lm_{3}^{j_3}\cdots \lm_{m}^{j_{m}}, \qquad f_3 =
\lm_{1}^{j_1-k}\lm_{2}^{j_2+k} \lm_{3}^{j_3} \cdots
\lm_{m}^{j_{m}},$$  are monomials and  $k \leq \min \{j_1, j_2\}$
is a non-negative  integer. Then $f_2$ is inessential for $f$.
\end{lemma}
\begin{proof}
Pick $\bfa = (a_1,\dots, a_m) \in \Trop^{(m)}$ and assume $a_1
> a_2$, then $f(\bfa) = f_1(\bfa) > f_2(\bfa), f_3(\bfa)$.
Conversely, if $a_2 > a_1$ then $f(\bfa) = f_3 (\bfa) > f_1(\bfa),
f_2(\bfa)$. When $a_1 = a_2$, $f(\bfa) = f_1(\bfa) \TrS f_3(\bfa)
= f_1(\bfa) \TrS f_2(\bfa) \TrS f_3(\bfa)$. Namely, $f_2$  is
inessential for $f$.
\end{proof}

\begin{remark}\label{rmk:Frobenius} A tropical polynomial $f \in
\Trop[\lm_1,\dots,\lm_m]$, written as a sum $\bigoplus_{i=1}^k
f_i$ of monomials $f_i$, satisfies the
 Frobenius property:
$$\bigg(\bigoplus_{i=1}^k f_i\bigg)^n = \bigoplus_{i=1}^k (f_i)^n$$ for any natural number
$n$, cf.~\cite[Corollary 3.23]{IzhakianRowen2007SuperTropical}. In
particular, $(a \TrS b)^n = a^n \TrS b^n$, for any $n \in \Net$.
\end{remark}


\subsection{Tropical matrices} It is standard that if $R$ is a semiring then we
can form the semiring $M_n(R)$ of all $n \times n$ matrices with
entries in $R$, where addition and multiplication are induced from
$R$. Accordingly, we define the semiring  of tropical matrices
$M_n(\Trop)$ over $\Trop = (\Trop, \TrS,  \TrP \ )$, whose unit is
the matrix
\begin{equation}\label{def:absUnitMat}  
\taI = \OP  \begin{array}{ccc}
  0 & \ldots & \tUniS \\
  \vdots & \ddots & \vdots \\
  \tUniS & \ldots & 0
\end{array} \CP \\
\end{equation}
and whose zero matrix is  $\taZ = (\tUniS) \taI$. We will consider
$M_n(\Trop)$ both as a semiring and also as a multiplicative
monoid. We denote tropical matrices as $\taA =(a_{ij})$ and use
the notation $a_{i j}$ for the entries of $\taA$. Clearly, since
$\Trop$ is commutative, $ \al \taA = \taA \al$, for any $\al \in
\Trop$ and $\taA \in M_n(\Trop)$.

\begin{note}\label{not:matGraph} In combinatorics,
an $n \times n $ tropical matrix $\taA \in M_n(\Trop)$ is used to
represent a weighted digraph $G = (E,V)$ with $n$ vertices $v_1,
\dots, v_n$, \cite{Lawler76}; we call this digraph the associated
digraph of $\taA$ and denote it by $G_\taA$. The edges $e \in E$
of $G_A$ are determined by pairs $(v_i,v_j)$ of vertices and the
weight $w(e)$ of an edge $e =  (v_i,v_j)$ of $G_A$ is the value of
the entry $a_{i j}$ in $A$. When $i = j$ the edge $e= (v_i,v_i)$
is called a self loop. Taking a power $\taA^i$ of $A$ is
equivalent to computing all the paths of length $i$ of maximal
weights on the associated graph $G_A$ of $\taA$,  \cite{brualdi}.
\end{note}

As usual, we define the \bfem{transpose} of $\taA = (a_{ij})$ to
be $\trn{\taA} = (a_{ji})$, and have the usual relation noted
here.
\begin{proposition}\label{prop:commutativeOfSimmMat}
$\trn{\OP \taA \taB \CP} =   \trn{\taB}   \trn{\taA} \ . $
\end{proposition}
\noindent  The proof follow easily from the commutativity and the
associativity of $ \TrS $ and $ \TrP $ over $\Trop$.

The \bfem{minor} $\taA_{ij}$ is obtained by deleting row $i$ and
column $j$ of $\taA$. We define the \bfem{tropical determinant} to
be
  \begin{equation}\label{def:tropicalDet}
 \per{\taA} = \bigoplus_{\sig \in S_n}  \OP  \aaa_{1\sig(1)}
 \cdots \aaa_{n\sig(n)} \CP,
\end{equation}
where $S_n$ is the set of all the permutations on $\{1,\dots,
  n\}$.   In terms of minors,  $\per{\taA}$ can be written equivalently as
$\per{\taA} = \bigoplus_{j}   \aaa_{i_o j} \per{\taA_{i_o j}},$
for some fixed index $i_o$. Indeed, in the classical sense, since
parity of indices are not involved, the tropical determinant is in
fact a permanent, which makes this definition purely
combinatorial.
\begin{remark}\label{rmk:powMat}
When $\taA$ has either a row (or a column) all of whose entries
are $\mInf$ then $\per{\taA} = \mInf$ and $\per{\taA^i} = \mInf$,
for each $i \in \Net$, since then $\taA^i$ also has either a row
(or a column) all whose entries are $\mInf$.
\end{remark}

A matrix $A \in M_n(\Trop)$ is said to be \textbf{tropically
singular} whenever the value of $|A|$, cf. Formula
\eqref{def:tropicalDet}, is attained by at least two different
permutations $\sig \in S_n$, otherwise $A$ is called
\textbf{tropically nonsingular}. The \bfem{adjoint} matrix
$\adj{\taA}$ of $\taA =(a_{ij})$ is defined as the matrix
$\trn{(a'_{ij})}$ where $a'_{ij}= \per{\taA_{ij}}$. When
$\per{\taA} \neq \mInf$ we use $\itaA$ to denote the tropical
quotient
$$\itaA \  :=  \ \frac{\adj{\taA}}{\per{\taA}} \ .$$
Note that the division in this definition is tropical division,
that is, substraction in the usual sense. The \bfem{multiplicative
trace} of $\taA$ is defined by the following formula.

\begin{equation}\label{eq:mtrace}
\mtrace{\taA} = \bigodot_i a_{ii} \ ,
\end{equation}
and therefore we always have $\per{A} \geq  \mtrace{\taA}$.

\begin{proposition}\label{prop:powOfMat} $\per{\taA^{n!}} \  = \  \mtrace{\taA^{n!}}$ for each  $\taA \in
M_n(\Trop)$.

\end{proposition}

\begin{proof} If $\per{\taA} = \mInf$  we are done, cf. Remark \ref{rmk:powMat}. Write
$B = (b_{ij})$ for $\taA^{n!}$  and assume that $\per{B}  >
\mtrace{B}$. In view of Note \ref{not:matGraph}, the graph $G_B$
of $B$ has a multicycle  $C$, of length $n$ whose weight is
greater than $\sum_{i} w(v_i,v_i)$. But since $B = \taA^{n!}$, the
weight of each self loop $(v',v')$ in $G_B$ is the maximal weight
over all paths of lengths $\ell \leq n$ from $v'$ to itself in
$G_A$ -- a contradiction.
\end{proof}

\begin{example}
Let us show by direct computation that $\per{\taA^{n!}} \  = \
\mtrace{\taA^{n!}}$ for the case of $n=2$.
$$
\begin{array}{llll}
  \Per{\vMat{a}{b}{c}{d}^{2}} & = &  \Per{\vMat{a^2 \TrS bc}{b(a \TrS
d)}{c(a \TrS d)}{d^2 \TrS bc}} & \\[1mm]
 &  = & (a^2 \TrS bc)(d^2 \TrS bc) \TrS
bc(a \TrS d)^2 & \\[1mm]
& =  & (a^2 d^2 \TrS  bc (a^2 + d^2) \TrS b^2 c^2  ) \TrS bc(a
\TrS d)^2 & \\[1mm]
& =  & a^2 d^2 \TrS  bc (a^2 + d^2) \TrS b^2 c^2   & =
\mtrace{\taA^{2}} \ .
\end{array}
$$
By the Frobenius property, cf. Remark \ref{rmk:Frobenius}, $bc(a
\TrS d)^2 = bc(a^2 \TrS d^2)$, so this component becomes
inessential and it is omitted.
\end{example}

\begin{remark} As a result of Proposition \ref{prop:powOfMat}, we can conclude
that not all tropical matrices have a square root; for example
take a matrix $\taA \in M_n(\Trop)$ with diagonal entries $=
\tUniS$ and all of its off-diagonal entries $\neq \tUniS$.
\end{remark}

A set $S$ of vectors $V_1, \dots, V_m \in \Trop^{(n)}$ is said to
be \bfem{linearly dependent} if there are $\al_1, \dots, \al_m \in
\Trop$, not all of them $\tUniS$, such that each coordinate of the
tropical sum $U = \bigoplus_t \al_t V_t$ is attained by at least
two different terms $\al_t V_{t}$, otherwise $S$ is linearly
independent. (In particular, if  $m >n $ then  $S$ is a  dependent
set \cite{zur05TropicalRank}.) The \bfem{rank}, $\rank{A}$,  of a
matrix $A \in M_n(\Trop)$  is the number of elements in a maximal
independent subset of rows.

\begin{theorem}[{\cite[Theorem 3.6]{zur05TropicalRank}}]\label{thm:IndependentToRegular}
An $n \times n$  matrix $A$ has rank $< n$ iff $A$ is tropically
singular.
\end{theorem}


 Using the theorem,  one
can easily check whether a matrix $A \in M_n(\Trop)$ has rank $n$
or not. For example, consider a matrix $A \in  M_2(\Trop) $  and
compute its tropical determinant:
$$\Per{\vMat{a}{b}{c}{d}}  = ad \TrS bc.$$
If $ad = bc$ then $A$ is tropically singular and thus has rank $<
2$, otherwise $A$ is nonsingular and is of rank $2$.


%

\subsection{Matrices of polynomials and polynomials of
matrices}\label{sec:MatricesofPoly}  We also look at the seimiring
and monoid of all matrices of polynomials. These are matrices
whose entries are elements of the polynomial semiring
$\Trop[\lm_1,\dots,\lm_m]$. We denote these matrices by
$M_n(\Trop[\lm_1,\dots,\lm_m])$ and for each $\bfa =
(a_1,\dots,a_m)\in \Trop^{(m)}$ we have the induced substitution
homomorphism:
\begin{equation}\label{eq:subsHomM}
\mat{\sHom_\bfa}:
 M_n(\Trop[\lm_1,\dots, \lm_m]) \ \To \  M_n(\Trop)
\end{equation}
given by sending $\pol{A} = (f_{ij}) \mapsto (f_{ij}(\bfa))$, cf.
Formula  \Ref{eq:subsHom}. We write $\pol{A}(\bfa)$ for the image
of $\pol{A}$ under the substitution homomorphism (to $\bfa \in
\Trop^{(m)}$). (Clearly, $M_n(\Trop) \hookrightarrow
M_n(\Trop[\lm_1,\dots, \lm_m])$ by sending $A = (a_{ij}) \mapsto
(f_{ij})$ where $f_{ij} = a_{ij}\lm_1^0 \cdots \lm_m^0 $ for all
$i,j = 1,\dots,n$.)

Inducing by  $e$-equivalence on $\Trop[\lm_1,\dots,\lm_m]$, we say
that two matrices $\pol{A} = (f_{ij})$ and $\pol{B} = (g_{ij})$ in
$M_n(\Trop[\lm_1,\dots,\lm_m])$ are \bfem{essentially equivalent},
denoted as
$$\pol{A}  \ \eqR \  \pol{B},$$ if
$f_{ij} \eqR g_{ij}$ for each $i$ and $j$. Therefore,  $\eqR$ is
an equivalence relation on $M_n(\Trop[\lm_1,\dots,\lm_m])$ for
which $\pol{A}$ and $\pol{B}$ are in a same class if and only if
$\pol{A}(\bfa) =  \pol{B}(\bfa) $ for each $\bfa \in \Trop^{(n)}$.

On the other hand, one can talk about polynomials, with
coefficients in $\Trop$, whose arguments are tropical matrices; we
denote this polynomial semiring by $\Trop[\Lm_1, \dots, \Lm_m]$
and write $F$, $G$, for its elements. Writing $\bfA$ for $(A_1,
\dots, A_m) \in M_n(\Trop)^{(m)}$, in the standard way, we define
the substitution homomorphism
$$\sHom_\bfA: \Trop[\Lm_1, \dots, \Lm_m] \ \To \  M_n(\Trop),$$
where $\sHom_\bfA : \bfA \mapsto F(\bfA)$. Viewing $\Lm_k =
({\lm_{ij}}^{(k)})$ as a matrix in $n ^2 $ indeterminates, we also
have the semiring homomorphism
$$\pmHom_{n,m} : \ \Trop[\Lm_1, \dots, \Lm_m] \ \To \  M_n(\Trop[\lm_1,\dots, \lm_M]), \quad M = m n^2,$$
given by sending $\pmHom_{n,m} :  \Lm_k \mapsto ({\lm_{ij}}^{(k)})
$. (Note that $\pmHom_{n,m}$ is not surjective.) Thus the
following diagram commutes.

\begin{equation}\label{diag:1}
{\setlength{\unitlength}{0.7mm}
\begin{picture}(60, 30)
  \put(-12,22){$\Trop[\Lm_1, \dots, \Lm_m]$}
  \put(56,22){$M_n(\Trop[\lm_1,\dots, \lm_M])$}
  \put(57,0){$ M_n(\Trop)$}

  \put(38,25){$\pmHom_{n,m}$}
  \put(40,12){$\sHom_\bfA$}
  \put(69,12){$\mat{\sHom_\bfa}$}

  \put(32, 23){\vector(1, 0){20}}
  \put(20, 20){\vector(2, -1){35}}
  \put(67, 20){\vector(0, -1){15}}
\end{picture}}
\end{equation}

\begin{remark}\label{rmk:matPoly}
Suppose $F,G \in \Trop[\Lm_1, \dots, \Lm_m]$ and consider their
images $\pmHom_{n,m}: F \mapsto \pol{A}$ and $\pmHom_{n,m}: G
\mapsto \pol{B}$.  In view of Diagram \Ref{diag:1}, when $\pol{A}
\eqR \pol{B}$ then $\pol{A}(\bfa) = \pol{B}(\bfa)$ for any $\bfa
\in \Trop^{(M)}$, and therefore $F(\bfA) = G(\bfA)$ for any $\bfA
= (A_1,\dots, A_m) \in M_n(\Trop)^{(m)}$.
\end{remark}

Remark \ref{rmk:matPoly} plays a key role in this paper and
provides the algebraic foundation for our semigroup applications,
especially in the study of semigroup identities.

\begin{remark} A non-essential polynomial $f$, that is, a polynomial $f$
that has  an inessential monomial, in $\Trop[\lm_1,\dots,\lm_m]$
can be essential as a polynomial of matrices; for example $\lm$ is
inessential monomial of $\lm^2 \TrS \lm \TrS 0$,  but the
corresponding polynomial $\Lm^2 \TrS \Lm \TrS I \in \Trop[\Lm]$ is
essential, that is $F = \Lm^2 \TrS \Lm \TrS I \neqR \Lm^2 \TrS I
=G $, since for $A = \vMat{\tUniS}{0}{0}{\tUniS}$,  $F(A) =
\vMat{0}{0}{0}{0}$while $G(A) = A$.

\end{remark}
\section{The Monoid $\mT$}

\subsection{Submonoids and subgroups of $M_n(\Trop)$} We begin by
describing some subsemigroups of $M_n(\Trop)$.

\begin{enumerate}\eroman
\item A tropical matrix with each row and column containing
exactly one entry $\neq \tUniS$ is called a permutation matrix.
The set of all tropical permutation matrices, which we denote by
$\tW_n$,  forms the \textbf{(affine) Weyl group}. It is the group
of units of $M_n(\Trop)$, that is the maximal subgroup  with
identity element $I$. $\tW_n$ contains the Abelian subgroup
$\tD_n$ of tropical diagonal matrices in $M_n(\Trop)$. \pSkip

\item The upper triangular matrices $U_n(\Trop)$ and the lower
triangular matrices $L_n(\Trop)$ are non-commutative  submonoids
of $M_n(\Trop)$ containing  $\tD_n$ as a submonoid. \pSkip

    \item
A  matrix, $\taA \in M_n(\Trop)$ is called \textbf{presymmetric}
matrix if it is symmetric about its anti-diagonal (i.e secondary
diagonal) and is said to be  \textbf{bisymmetric} if it is both
symmetric and presymmetric. It is easy to see that the bisymmetric
$2 \times 2$ matrices form a commutative submonoid of $\mT$.
 \pSkip

\item A tropical matrix $A \in \mT$ with all diagonal entries
equal $0$ and off-diagonal entries $\leq 0$ is an idempotent
matrix, i.e. $A^2 = A$. It is easy to see that the product of any
two such matrices is equal to their sum and thus commute. Thus the
collection of all such matrices is a submonoid of $\mT$ denoted by
$\tN_2$. This is a commutative monoid and $\per{A} = \mtrace{A}
=0$ for each $A \in \tN_2$.
\end{enumerate}


\subsection{Von-Neumann  regularity} Regularity of semigroup and invertibility of
their elements have several different notions,  in this paper we
use  Von-Neumann's notion.

\begin{definition} Let $(\tS, \cdot \ )$ be a semigroup, an
element $y \in \tS$ is called a \textbf{generalized inverse} of $x
\in \tS$ if
$$ x \, y \, x  \ = \ x \quad \text{and} \quad y \, x \, y \  = \ y \ , $$
in the case that $x$ has an inverse we say that $x$ is regular in
$\tS$. A semigroup $\tS$ is said to be a \textbf{regular
semigroup}, in the Von-Neumann sense, if every element $x \in \tS$
has at least one generalized inverse $y \in \tS$.
\end{definition}

\begin{lemma}\label{lem:inver} Each matrix $A \in \mT$ with $\per{A} = \tUniS$
has a generalized inverse.
\end{lemma}

\begin{proof} Since $\per{A} = \tUniS$, then $A$ has either a row or a column whose entries  are all
$\tUniS$, so it has one of the following forms:
\begin{equation*}\label{eq:inv.1}
    A_1 = \vMat{a}{b}{\tUniS}{\tUniS}, \quad A_2 =
    \vMat{\tUniS}{\tUniS}{a}{b}, \quad
    A_3 = \vMat{a}{\tUniS}{b}{\tUniS}, \quad A_4 =
    \vMat{\tUniS}{a}{\tUniS}{b},
\end{equation*}
where $a,b \in \Trop$. Correspondingly, we specify their inverses
to be:
\begin{equation*}\label{eq:inv.1}
    B_1 = \vMat{-a}{\tUniS}{-b}{\tUniS}, \quad B_2 =
    \vMat{\tUniS}{-a}{\tUniS}{-b}, \quad
    B_3 = \vMat{-a}{-b}{\tUniS}{\tUniS}, \quad B_4 =
    \vMat{\tUniS}{\tUniS}{-a}{-b} \ .
\end{equation*}
Note that when $a$ or $b$ is $\tUniS$, then, respectively, the
terms $-a$ or $-b$ are replaced by $\tUniS$.
\end{proof}

\begin{theorem}\label{thm:V.N.for2by2rerular} Assume $\per{\taA} \neq
\mInf$, then $\itaA$ is a generalized inverse of $\taA \in \mT$.
\end{theorem}
\begin{proof}
Suppose  $\taA = \vMat{a}{b}{c}{d}$ with $\per{\taA} \neq \mInf$,
 then $\itaA = \vMat{d}{b}{c}{a}/ \per{A}$. Computing the product
$\taA \itaA$
 $$\taA \itaA = \vMat{a}{b}{c}{d} \vMat{d}{b}{c}{a}/ |A| =
    \vMat{0}{\frac{ab}{|\taA|}}{\frac{cd}{|\taA|}}{0} $$
we get
 $$\taA \itaA \taA = 
     \vMat{0}{\frac{ab}{|\taA|}}{\frac{cd}{|\taA|}}{0} \vMat{a}{b}{c}{d} =
     \vMat{a(0 \TrS \frac{bc}{|\taA|} )}{b(0 \TrS \frac{ad}{|\taA|})}
     {c(0 \TrS \frac{ad}{|\taA|})}{d(0 \TrS \frac{bc}{|\taA|})}.$$
The proof is then derived from the relations  $bc, ad \leq
\per{\taA}$. The relation  $\itaA \taA \itaA = \itaA$ is proved in
the same way.
\end{proof}

\begin{corollary}\label{cor:m2reg}
$\mT$ is a regular semigroup.
\end{corollary}
\begin{proof}
Use Lemma \ref{lem:inver} and Theorem \ref{thm:V.N.for2by2rerular}
according to when $\per{\taA} = \mInf$ or not.
 \end{proof}

\begin{remark}

The relations  $\taA \itaA \taA = \taA$ and $\itaA \taA \itaA
 = \itaA$ fail  for $\taA \in M_n(\Trop)$. For example take  $\taA \in
 M_3(\Trop)$ to be
$$\taA = \vvMat
 {-4}{4}{-2}
 {0}{-1}{-3}
 {1}{-2}{-3}
\quad \text{then ,} \quad  \per{\taA} = 2 \quad \text{and} \quad
\itaA = \vvMat {-6}{-1}{-1} {-4}{-3}{-4} {-2}{3}{2} \ ,
$$
the product is then
$$
\taA \itaA \taA = \vvMat{1}{4}{-2}{0}{-1}{-3}{ 1}{ -1}{-3} \  \neq
\ \taA \ .
$$

The monoid $M_n(\Trop)$, $n > 2$, is not regular.  Indeed, let
$\bool$ denote the 2 element boolean semiring. Then there is a
semiring homomorphism
$$ \psi : \Trop \To \bool,$$
given by sending $-\infty \mapsto 0$ and $a \mapsto 1$ for any $a
\in \Real$, which induces a surjective monoid homomorphism
$$\mat{\psi}: M_n(\Trop) \To M_n(\bool). $$
But it is known that $M_n(\bool)$ in not regular if $n \geq 3$
(see \cite[Chapter 2]{kim82}) and so neither is $M_n(\Trop)$ since
regularity is preserved by surjective morphisms.
\end{remark}

\begin{proposition} When $\per{\taA} \neq \tUniS$,
$\taA \itaA = \itaA \taA$ if and only if $\taA$ is presymmetric,
for any $\taA \in \mT.$
\end{proposition}
\begin{proof} Computing the products
 $$\taA \itaA =
    \vMat{0}{\frac{ab}{|\taA|}}{\frac{cd}{|\taA|}}{0}  \quad \text{ and
    } \quad  \itaA \taA  =
    \vMat{0}{\frac{bd}{|\taA|}}{\frac{ac}{|\taA|}}{0} \ ,  $$  as far as $\per{\taA}
\neq \tUniS$, one can see that  $\taA \itaA = \itaA \taA$ if and
only if $a =d$.
\end{proof}

\section{Semigroup Identities On $\mT$}

Our main result in this section is that the monoid $\mT$ admits a
non-trivial semigroup identity. We also show that there are
 semigroup identities for other submonoids of  $\mT$ like
for triangular tropical matrices.

\subsection{Semigroup identities}
Assuming  $(\tS, \cdot \ )$ is a semigroup with an identity
element $1$, we write $x^i$ for the $x \cdot x \cdots x$ repeated
$i$ times and identify $x^0$ with $1$.

Let $X$ be a countably infinite set of ``variables". A semigroup
identity is a formal equality $u=v$ where $u$ and $v$ are in the
free semigroup $X^{+}$ generated by $X$. For a monoid identity, we
allow $u$ and $v$ to be the empty word as well. A semigroup $S$
satisfies the semigroup identity $u=v$ if for every morphism
$f:X^{+} \rightarrow S$, one has $uf=vf$. Let $I$ be a set of
identities. The set of all semigroups satisfying every identity in
$I$ is denoted by $V[I]$ and is called the variety of semigroups
defined by $I$. It is easy to see that $V[I]$ is closed under
subsemigroups, homomorphic images and direct products of its
members. The famous Theorem of Birkhoff says that conversely, any
class of semigroups closed under these three operations is of the
form $V[I]$ for some set of identities $I$.

\begin{example}




     The tropical Weyl group $\tW_n \subset  M_n(\Trop)$
    satisfies the  identity  $A^{n!}B^{n!} = B^{n!}A^{n!}$, since
    $A^{n!} \in \tD_n$ for each $A \in \tW_n$ (cf. Proposition
    \ref{prop:powOfMat}), which is an Abelian group.

\begin{remark}

The preceding example shows that the group of units of
$M_n(\Trop)$ is a virtually Abelian group.

%

\end{remark}

\end{example}



\begin{remark}\label{rmk:idToMonom} In view of Subsection \ref{sec:MatricesofPoly},
for the case of $\tS = M_n(\Trop)$ we can identify any semigroup
identity $u=v$ with a pair of monomials $ H^{u}, H^{v} \in \
M_n(\Trop)[\Lm_1, \Lm_2]$ whose powers are determined by the words
$u$ and $v$ and their coefficients are in $M_n(\Trop)$.
\end{remark}

\subsection{The submonoid of triangular matrices }
In this section we prove that
\begin{equation}\label{eq:biId.2}
AB^2A \ AB  \ AB^2A = AB^2A \  B A \ AB^2A \ ,
\end{equation}
is a semigroup identity for the submonoid  $\uT \subset \mT$ of
upper triangular tropical matrices.

\begin{remark}\label{rmk:rank1} In the case when $A \in \mT$ is of rank $1$ it is easy to verify
that $A^2 = \al A$, for some $\al \in \Trop$; therefore, the
semigroup identity  \Ref{eq:biId.2} is satisfied in $\mT$ whenever
$A$ or $B$ is of rank $1$. (To see that, just extract the scalar
multiplier to obtain the equality.) By the same argument, when
$\pol{A}$ or $\pol{B}$ is of rank $1$, they also satisfy relation
\Ref{eq:biId.2}.
\end{remark}

\begin{lemma}\label{lem:zeroTR}
Suppose  $\pol{A}, \pol{B} \in \upT$ are of the form
\begin{equation}\label{eq:zeroTR}
\pol{A} = \vMat{0}{\lm_1}{\tUniS}{\lm_2}, \qquad
  \pol{B} = \vMat{0}{\lm_3}{\tUniS}{\lm_4}
\end{equation}
then $u(\pA,\pB) \eqR v(\pA,\pB)$, where $u=AB^2A \ AB  \ AB^2A$
and $v=AB^2A \ BA  \ AB^2A$.
\end{lemma}

\begin{proof}

Compute the products $\pF= u(\pA,\pB)$ and $\pG = v(\pA,\pB)$,
write $\pF= (f_{ij})$ and $\pG= (g_{ij})$, and consider the
$f_{ij}$'s and the $g_{ij}$'s as flats polynomials in 4
indeterminates, $\lm_1, \dots , \lm_4$. It is easy to verify that
\pSkip
$$f_{11} = g_{11} = 0, \quad
f_{21} = g_{21} = \tUniS ,  \quad \text{ and}  \quad
f_{22} = g_{22} = \lm_2^5 \lm_4^5 .$$

Writing $f_{12} = h_{12}  \TrS  \al_{12}$ and $g_{12} = h_{12}
\TrS \bt_{12}$,  where by direct computation we have,
$$h_{12} = \lm _1+\lm _2 \lm _3+\lm _2 \lm _3 \lm _4+\lm _1 \lm _2 \lm _4^2+\lm _1 \lm _2^3 \lm _4^3+\lm _2^4 \lm _3 \lm _4^3+\lm _2^4 \lm _3 \lm _4^4+\lm _1 \lm _2^4 \lm _4^5 \ , $$
$$\al_{12} = \lm _2^2 \lm _3 \lm _4^2+\lm _1 \lm _2^2 \lm _4^3
 \ , \qquad  \text{and} \qquad \bt_{12} = \lm _1 \lm _2^2 \lm _4^2+\lm _2^3 \lm _3 \lm _4^2 \ . $$
The diagrams below, together with Lemma \ref{lem:essential}, show
that all the terms of $\al_{12}$ and $\bt_{12}$ are inessential
for $h_{12}$, and thereby  $f_{12} \eqR g_{12}$: \noindent
$$  \begin{array}{l}    \al_{12} = \end{array}
\begin{array}{ccc}
        \lm _2 \lm _3 \lm _4 && \lm _1 \lm _2 \lm _4^2  \\
         | & &  | \\
  \lm _2^2 \lm _3 \lm _4^2 &  \TrS  & \lm _1 \lm _2^2 \lm _4^3  \\
   | & &  | \\
 \lm _2^4 \lm _3 \lm _4^4&& \lm _1 \lm _2^4 \lm _4^5
 \\
\end{array}
\qquad
 \begin{array}{l}   \bt_{12} = \end{array}
\begin{array}{ccc}
       \lm_1 && \lm _2 \lm _3  \\
         | & &  | \\
 \lm _1 \lm _2^2 \lm _4^2 &  \TrS  &\lm _2^3 \lm _3 \lm _4^2 \\
   | & &  | \\
  \lm _1 \lm _2^3 \lm _4^3  && \lm _2^4 \lm _3 \lm _4^3
 \\
\end{array} \ .
$$
The upper and the lower rows specify the monomials in $h_{12}$
that make respectively each terms of $\al_{12}$ and $\bt_{12}$ to
be inessential.  Taking all together, $f_{ij} \eqR g_{ij}$ for all
$i,j = 1,2$, and thus $\pol{F} \eqR \pol{G}$.
\end{proof}

\begin{theorem}\label{thm:idTr} The submonoid $\uT$ of upper triangular
tropical matrices admits the semigroup identity $$ AB^2A \ AB  \
AB^2A = AB^2A \  B A \ AB^2A \ .$$
\end{theorem}
\begin{proof} Take $A, B \in \uT$, if one of them is of rank $1$
we are done by Remark \ref{rmk:rank1}; otherwise, we can divide
$A$ by $a_{11}$ and $B$ by $b_{11}$ to have matrices of the form
\Ref{eq:zeroTR}, then in the view of Remark \ref{rmk:idToMonom}
the proof is completed by Lemma \ref{lem:zeroTR}.
\end{proof}

\begin{corollary} The submonoid of lower triangular matrices $\lT$ admits the semigroup identity \Ref{eq:biId.2}.
\end{corollary}
\begin{proof}
Immediate by Theorem \ref{thm:idTr} and the fact that $\lT$ is
conjugate to the monoid $\uT$.
\end{proof}

\subsection{A semigroup identity on $\mT$ }

\begin{lemma}\label{lem:zeroM}
Suppose  $\pol{A}, \pol{B} \in \mpT$ are of the form
\begin{equation}\label{eq:zeroM}
\pol{A} = \vMat{0}{\lm_1}{\lm_2}{\lm_3}, \qquad
  \pol{B} = \vMat{0}{\lm_4}{\lm_5}{\lm_6}
\end{equation}
satisfying the restriction   $\per{A} = \lm_3$ and $\per{B} =
\lm_6$, then $u(\pA,\pB) \eqR v(\pA,\pB)$, where
$u=AB^{2}AABAB^{2}A$ and $v=AB^{2}ABAAB^{2}A$.
\end{lemma}
\noindent The proof of the lemma is long and technical and is
proved in detail as Lemma \ref{lem:zeroMA} in the Appendix.

\begin{theorem}\label{thm:globalId} The monoid $\mT$ admits the semigroup identity
\begin{equation}\label{eq:
Id.2} A^2B^4A^2 \ A^2B^2  \ A^2B^4A^2 =A^2B^4A^2 \  B^2 A^2 \
A^2B^4A^2 \ .
\end{equation}
\end{theorem}

\begin{proof}
Suppose $A, B \in \mT$. If one of these matrices is of rank $1$ we
are done by Remark \ref{rmk:rank1}. Otherwise, taking their
squares, by Proposition \ref{prop:powOfMat}, $\per{A^2} =
\mtrace{A^2}$ and $\per{B^2} = \mtrace{B^2}$, so we can divide
their squares to have matrices of the form $\Ref{eq:zeroM}$.
Consider the matrices over $M_2(\Trop[\lm_1,\dots,\lm_6])$
corresponding to these squares, i.e. matrices of the form
$\Ref{eq:zeroM}$ whose entries are monomials. In view of Remark
\ref{rmk:idToMonom}, the proof  is then completed by Lemma
\ref{lem:zeroM}.
\end{proof}

\section{The bicyclic monoid embeds in $\mT$}



It is well known that the monoid of matrices $M_{n}(K)$ over a
field $K$ is semisimple \cite{CP} and in particular does not have
a copy of the bicyclic monoid $\tB$ as a subsemigroup. See
\cite{Okninski}, Chapter 2 for the basic semigroup structure of
$M_{n}(K)$. In contrast to this, we prove in this section that
$\tB$ has a faithful representation in $\mT$. This explains in
part, the identity in Theorem \ref{eq:biId.2}. This is Adjan's
identity for the bicyclic monoid and is the shortest identity
satisfied by $\tB$, cf. \cite{Adjan}.

Although in this paper we do not study properties of semigroups by
their actions on tropical spaces,  we open by presenting the
tropical analogue of a linear representation.
\subsection{Tropical Linear Representations}

Considering $\Trop^{(n)}$ as a space of finite dimension we denote
the (tropical) associative semialgebra of all tropical linear
operators on $\Trop^{(n)}$ by $L(\Trop^{(n)})$. These linear
operators can be represented as matrices (in some basis) and this
establishes an isomorphism between $L(\Trop^{(n)})$ and the matrix
algebra $M_n(\Trop)$. Therefore, we can identify $L(\Trop^{(n)})$
with $M_n(\Trop)$. Recall that our ground structure is a semiring
and thus, the notions of spaces, operators, and algebra are the
corresponding notions \cite{Lallement}.

 A finite dimensional \textbf{tropical linear representation} of a
semigroup $\tM$, over $\Trop^{(n)}$, is a semigroup homomorphism
$$\rep: \  \tM \To L(\Trop^n)$$
%
%
(The space $\Trop^{(n)}$ can be replaced by other tropical spaces,
but to clarify the exposition we focus on $\Trop^{(n)}$.) When
$\rep$ is a one-to-one homomorphism, then the representation is
called \textbf{faithful}. As in classical representation theory
one should think of a representation as a tropical linear
\textbf{action} of $\tM$ on $\Trop^{(n)}$ (since to every $a \in
\tM$, there is associated a tropical linear operator $\rep(a)$
which acts on $\Trop^{(n)}$).

\subsection{A tropical representation of the bicyclic monoid  }

 The monoid, $\tB = \langle a, b \rangle$, generated by two elements $a$ and
 $b$  satisfying the one relation
\begin{equation}\label{eq:biRole}
ab = 1,
\end{equation} where $1$ is the identity element, is
called the \textbf{bicyclic monoid}. The elements $x,y \in \tB $
of $\tB$ are called words (or strings) over $a$ and $b$. It is
well known that every element of $\tB$ is equal to a unique word
of the form $x = b^i a^j, \qquad i,j \in \Int_+ $ \cite{CP}. As
usual, we identify the elements $a^0$ and $b^0$ with the identity
element $1$ of $\tB$.

We start by recalling another representation of the elements of
$\tB$ which helps us later to formulate a faithful tropical linear
representation of $\tB$.

%
%

Let $\tS$ denote $\mathbb{N} \times \mathbb{N}$. We define a
binary operation $\ast$ on $\tS$ by the following formula.

\begin{equation}\label{eq:bi2powers}
* : ((i,j),(h,k)) = \left\{ \begin{array}{lcc}
                (i+h-j,k), &  & j \leq h, \\[1mm]
                (i,j-h+k,)  &  & j > h, \\
                              \end{array}
       \right.
\end{equation}

The following proposition is classical \cite{CP}.

\begin{proposition}\label{thm:isoOfBiCyclic.1} Given  a bicyclic monoid $\tB$,
the map $\phi: \tB \to \tS$, where $\phi: b^i a^j \mapsto (i,j)$,
is  a monoid isomorphism.
\end{proposition}

\begin{proof}
Clearly, $\phi$ is  a bijective, and
$$ \phi(( b^i a^j)( b^h a^k)) =
\left\{
\begin{array}{l}
    \phi(b^{i+h-j}a^k) \\[1mm]
    \phi(b^{i}a^{k+j-h})  \\
\end{array}
\right.  = \ \left\{
\begin{array}{ll}
    (i+h-j,k), & h  \geq j ; \\[1mm]
    (i,{k+j-h}), & h < j; \\
\end{array}
\right.  = (\bar{i},j) * (\bar{h},k) \ .$$
\end{proof}

We use this isomorphism to define a tropical linear representation
of $\tB$. Let $\tU_2$  be the subsemigroup of $U_2(\Trop)$, the
monoid of $2 \times 2$ upper triangular tropical matrices,
generated by the two elements
\begin{equation}\label{eq:genTropicalBiCyclic}
A  = \vMat{1^{-1}}{1}{\tUniS}{1} \qquad  \text{ and } \qquad
       B = \vMat{1}{1}{\tUniS}{1^{-1}},
\end{equation}
we write $1^{-1} = \frac{0}{1}$, which is just $-1$ in the usual
sense. Having these generators for $\tU_2$, $A^j$ and $B^i$, for
$i,j \in \Net$, can written as
$$  \ A^j =
\vMat{j^{-1}}{j}{\tUniS}{j},  \qquad  B^i =
\vMat{i}{i}{\tUniS}{i^{-1}}, $$ and they satisfy the following
relations:
\begin{equation}\label{eq:tropicalBiCyclic}
 E = A B =
\vMat{0}{0}{\tUniS}{0}, \qquad BA = \vMat{0}{2}{\tUniS}{0}, \qquad
B^i A^j = \vMat{j^{-1}i}{ij}{\tUniS}{j i^{-1}}  \ ,
\end{equation}
where   $i^{-1}j = \frac{j}{i}$ .

\begin{corollary}\label{thm:tropicalBiCyclic.2}  $\tU_2$ is a bicyclic
monoid.
\end{corollary}
\begin{proof}
Immediate from \Ref{eq:tropicalBiCyclic} and the fact that $E$ is
the identity element of $\tU_2$.
\end{proof}
Denoting the tropical semiring, having the addition $\TrS$ and
multiplication $\TrP$,  over $\Int \cup \{ \tUniS \}$ as $\eInt$,
the operation of the monoid $\tS$, cf. Eq.~\Ref{eq:bi2powers}, is
translated naturally to product of matrices in over $\tU_2 \subset
U_{2}(\eInt)$.

\begin{proposition}\label{thm:isoOfBiCyclic.2}
The map
$$\psi \ : \ \tS \To \tU_2, \qquad \psi: (\bar{i},j) \mTo
\vMat{j^{-1}i}{ij}{\tUniS}{j i^{-1}},$$ is a monoid isomorphism.
\end{proposition}
\begin{proof}
Assume $\vMat{a}{b}{\tUniS}{c} = \vMat{j^{-1}i}{ij}{\tUniS}{j
i^{-1}}, $
then $j = a^{-1}i$, $j =ci$, and $b = ij$. Accordingly, $a =
c^{-1}$ and thus $\psi$ is bijective.
\end{proof}

Finally, we define the monoid isomorphism
\begin{equation*}\label{eq:repOfBiCyclic} \rep \ : \ \tB \To
\tU_2 \
\end{equation*}
as the composition $\rep = \psi \circ \phi$ to get:
 \begin{theorem}\label{thm:repOfBiCyclic}
$\rep$ is  a faithful linear representation of $\tB$.
\end{theorem}

\begin{proof}
$\rep$ is a composition of monoid isomorphisms onto a matrix
monoid, and thus is a faithful tropical linear representation.
\end{proof}

\begin{corollary}\label{cor:bicyclic} The bicyclic monoid $\tB$ satisfies the semigroup identity
\Ref{eq:biId.2}, i.e.
$$  x y^{2} x \; x y \; x y ^{2}x =x y ^{2}x \; y x \; xy ^{2}x   $$
for any $x,y \in \tB$.
\end{corollary}
\begin{proof} Immediate by Corollary \ref{cor:m2reg} and Theorem
\ref{thm:repOfBiCyclic}.
\end{proof}

As mentioned above the semigroup identity in Corollary
\ref{cor:bicyclic} is known as Adjan's identity for the bicyclic
monoid \cite{Adjan}. In this paper we have provided an alternative
approach  for proving this semigroup identity and maybe other
semigroup identities. Tropical representation theory thus is
useful for studying properties of semigroups where classical
representation theory does not have anything to say. A geometric
point of view for the identities of a bicyclic monoids is provided
in \cite{pastijnBicyclic}.

Note that the morphism in Corollary \ref{cor:bicyclic} is not a
monoid morphism- it does not take the identity element of $\tB$ to
the identity element of $\mT$. In fact, no such faithful monoid
morphism exists. It is not difficult to see that the
$\mathcal{D}$-class of 1 in $\mT$ is the Weyl group $\tW_n$. Since
$\tB$ is a bisimple monoid that is not a group, no faithful monoid
morphism exists. See \cite{TropGreen} for more information on
Green's relations in full tropical matrix monoids.

\subsection{Remarks and Open Problems}

We have proved that $\mT$ satisfies a non-trivial semigroup
identity and that begs the question about whether $M_{n}(T)$
satisfies non-trivial identities for all $n>2$. We conjecture that
this is so, but have not been able to prove this as of yet.

For $n=2$, the connection between the monoid of upper triangular
$2 \times 2$ tropical matrices and the bicyclic monoid is deeper
than we've indicated in this paper. In fact the monoid of upper
triangular matrices of rank 2 form an inverse monoid that is
isomorphic to the monoid of partial shifts of the real line, just
like $\tB$ is isomorphic to the monoid of partial shifts on the
natural numbers \cite{TropGreen}. It was this connection that lead
us to try Adjan's identity on the submonoid of upper triangular
matrices in $\mT$ and eventually to the identity for all of $\mT$
that we found in this paper. We would like to clarify the exact
relationship between the bicyclic monoid and the monoid of upper
triangular $2 \times 2$ matrices further. We ask if they generate
the same variety, that is, if they satisfy exactly the same
identities.

The structure of upper triangular full rank $n \times n$ tropical
matrices is illimunated in \cite{TropGreen}. It is a block group,
that is a monoid in which each $\mathcal{R}$ and $\mathcal{L}$
class have at most one idempotent, but is not an inverse monoid if
$n>2$. This has made finding an identity difficult computationally
for this monoid. Passing to the monoid of all $n \times n$
matrices is also difficult.

Another reason for conjecturing that $M_{n}(\Trop)$ satisfies a
non-trivial identity for all $n$ is that every finite subsemigroup
of $M_n(\Trop)$ has polynomial growth \cite{Gaub, Simon}. In
particular, the free semigroup on 2 generators is not isomorphic
to a subsemigroup of $M_n(\Trop)$. While Shneerson \cite{Shn} has
given examples of polynomial growth semigroups that do not satisfy
any non-trivial identity (no such example exists for groups by
Gromov's Theorem \cite{Grom}, we feel that this is not the case
for $M_{n}(\Trop)$.
%

\bibliographystyle{abbrv}


\begin{thebibliography}{1}

\bibitem{Adjan}
S.~I. Adjan,
\newblock {\em Defining relations and algorithmic problems for groups and
  semigroups}.
\newblock Number~85. Proceeding of the Steklov Institute of Mathematics,
  American Mathematical Society, 1967.

\bibitem{brualdi}
R.~A. Brualdi and H.~J. Ryser,
\newblock {\em Combinatorial matrix theory}.
\newblock Cambridge University Press, 1991.

\bibitem{CP}
A.H.~Clifford and G.B.~Preston,
\newblock {\em The Algebraic Theory of Semigroups},
\newblock AMS,Providence, R.I., , Volume 1, 1961, Volume 2, 1967.

\bibitem{Gaub}
S.~Gaubert and R.D.~Katz,
\newblock {Reachability problems for products of matrices in
semirings},
\newblock {\em Int. J. of Alg. and Comp.}, Vol. 16 no. 3(2006), 603-627.

\bibitem{Grom}
\newblock {\em Groups of polynomial growth and expanding maps}
\newblock Publ. Math IHES 53 (1981), 53-73.

\bibitem{IMS}
I.~Itenberg, G.~Mikhalkin, and E.~Shustin.
\newblock {\em Tropical algebraic geometry}, volume~35.
\newblock Birkhauser, 2007.
\newblock Oberwolfach seminars.

\bibitem{zur05TropicalRank}
Z.~Izhakian.
\newblock The tropical rank of a tropical matrix.
\newblock Preprint at arXiv:math.AC/0604208, 2005.

\bibitem{TropGreen}
Z.~Izhakian and  S.W. Margolis
\newblock {\em Green's Relations on the monoid of all tropical
matrices}.
\newblock To appear.

\bibitem{IzhakianRowen2007SuperTropical}
Z.~Izhakian and L.~Rowen,
\newblock Supertropical algebra,
\newblock preprint at arXiv:0806.1175, 2007.

\bibitem{kim82}
K.~H. Kim.
\newblock {\em Boolean Matrix Theory and Applications}, volume~70 of {\em
  Monographs and Textbooks in Pure and Applied}.
\newblock Marcel Dekker, New York, 1982.

\bibitem{Lallement}
G.~Lallement.
\newblock {\em Semigroups and Combinatorial Applications}.
\newblock John Wiley \& Sons, Inc., New York, NY, USA, 1979.

\bibitem{Lawler76}
E.~L. Lawler.
\newblock {\em Combinatorial Optimization: Networks and Matroids}.
\newblock Holt, Rinehart, and Winston, 1976.

\bibitem{Okninski}
J.Okninski,
\newblock {\em Semigroups of Matrices}
\newblock World Scientific, Singapore, 1998.

\bibitem{pastijnBicyclic}
F.~J. Pastijn.
\newblock Polyhedral convex cones and the equational theory of the bicyclic
  semigroup.
\newblock {\em J. Austra. Math. Soc}, 81:63--96, 2006.

\bibitem{pin98}
J.-E. Pin.
\newblock Tropical semirings.
\newblock {\em Cambridge Univ. Press, Cambridge}, 11:50--69, 1998.
\newblock Publ. Neton Inst. 11, Cambridge Univ.

\bibitem{RST}
J.~Richter-Gebert, B.~Sturmfels, and T.~Theobald.
\newblock First steps in tropical geometry.
\newblock {\em Idempotent mathematics and mathematical physics}, pages
  289--317, 2005.
\newblock Contemp. Math., Amer. Math. Soc., Providence, RI, 377.

\bibitem{Shn}
L.~Shneerson
\newblock{\em Identities in finitely generated semigroups of polynomial growth}
\newblock J. Algebra  154  (1993),  no. 1, 67--85.

\bibitem{Simon}
I.~Simon,
\newblock {\em Recognizable sets with multiplicities in the tropical
semiring}, in MFCS 88, editors, M. Chytil, L. Janiga, V. Koubek,
Lecture Notes in Computer Science, Number 324, Springer, 107-120,
1988.


\end{thebibliography}

\vskip 1cm
\section*{Appendix A}

To clarify the exposition, instead of $\lm_1, \dots, \lm_6$, we
use the letters $a,b,c,x,y,$ and $z$, to denote the variables of
matrices of polynomials $\pol{A}, \pol{B} \in M_2(\Trop[\lm_1,
\dots, \lm_6] )$.
\begin{lemma}\label{lem:zeroMA} Given two matrices $\pol{A},\pol{B} \in M_2(\Trop[a,b,c,x,y,z])$ of the form
 \begin{equation*}\label{eq:app1} \pol{A} = \vMat{0}{a}{b}{c}
\qquad \text{ and } \qquad
    \pol{B} = \vMat{0}{x}{y}{z},
\end{equation*} assuming $\per{\pol{A}} = c$ and $\per{\pol{B}} =
    z$, then  \begin{equation}\label{eq:app1}  \pol{A} \pol{B}^2\pol{A}^2 \pol{B}
\pol{A} \pol{B}^2\pol{A} \ \eqR \ \pol{A} \pol{B}^2\pol{A} \pol{B}
\pol{ A}^2 \pol{B}^2\pol{A} \ .\end{equation}
\end{lemma}
By writing $\per{A} = c$ and $\per{B} = z$ we actually mean that
for any substitution of $\bfa \in \Trop^{(6)}$ we have
$\per{\pol{A}}(\bfa) = c(\bfa)$ and $\per{\pol{B}}(\bfa) =
    z(\bfa)$ .

\begin{proof} Taking the two products $\pol{F} = \seq^{\l}(\pA,\pB)$
 and $\pol{G} = \seq^{\r}(\pA,\pB)$  respectively for  the left and the right hand side of
 \eqref{eq:app1},
    written respectively as $\pol{F} = (f_{ij})$ and  $\pol{G} = (g_{ij})$,
and considering each $f_{ij}$ and $g_{ij}$, $i,j = 1,2$,  as
polynomials in $\Trop[a, b, c, x, y, z]$ we show that each pair
$f_{ij}$ and $g_{ij}$ of polynomials are $e$-equivalent.  (In fact
 $ \pol{F}$ and $\pol{G}$ are just elements  of $M_2(\Trop[a,b,c,x,y,z]$ whose preimages in
 $\Trop[\Lm_1, \Lm_2]$ under $\pmHom_{n,m}$ are monomials,
 cf. Diagram \eqref{diag:1}.)
%
%
%
%
 Note that by this construction all the coefficients of $f_{ij}$
and $g_{ij}$ are constantly $0$, thus all are flat polynomials.

In this view we prove that entry-wise
$$ \pol{F} \ \eqR \ \pol{G} \ . $$
 To do so, for each pair $f_{ij}$ and $g_{ij}$ we write $f_{ij} = h_{ij}
+ \al_{ij}$ and $g_{ij} = h_{ij} + \bt_{ij}$, where $h_{ij}$,
$\al_{ij}$,  and $\bt_{ij}$ are  polynomials in $\Trop[a, b, c, x,
y, z]$, and show that each monomial in $\al_{ij}$  and $\bt_{ij}$
is inessential with respect to~$h_{ij}$.

Recall that by the hypothesis of the lemma,
\begin{equation}\label{eq:a1}
    \per{\pol{A}} = c \geq ab \quad \text{ and } \quad \per{\pol{B}} = z \geq
    xy \ .
\end{equation}
Using this property, whenever a monomial $f_s$ of $f$ is
``greater" that $f_t$, that is $f_s(\bfa) \geq f_t(\bfa)$ for each
$\bfa \in \Trop^{(6)}$, then $f_t$ is inessential. We call these
``greater" monomials \bfem{dominant monomials} and mark them in
the  text using the bold font.

We complete the proof by observing the different entries of the
matrices $\pol{F}$ and $\pol{G}$ and identifying for each case the
monomials which make $\al_{ij}$ and $\bt_{ij}$ inessential for
$h_{ij}$, respectively. These monomials are labeled in $h_{ij}$ by
$\underline{\ \square\ }^{(\al.\#no)}$ and $\overline{\ \square \
}^{(\bt.\#no)}$ according to their numbering in $\al_{ij}$ and
$\bt_{ij}$ which appear in the diagrams below.

To make the next technical exposition clearer, the monomials in
the  equations below are ordered lexicographically.
\begin{description}
    \item[$\underline{i=1, j =1}$] $ $ \pSkip

\noindent $h_{11} =1 \TrS a b \TrS b x \TrS a b^2 x \TrS b c x
\TrS b^2 x^2 \TrS a b^3 x^2 \TrS b^2 c x^2 \TrS b^3 x^3 \TrS b^3 c
\ x^3 \TrS a y \TrS a^2 b y \TrS a c y \TrS x y \TrS a b x y \TrS
a^2 b^2 x y \TrS c x y \TrS a b c x y \TrS c^2 x \ y \TrS b x^2 y
\TrS a b^2 x^2 y \TrS a^2 b^3 x^2 y \TrS b c x^2 y \TrS a b^2 c
x^2 y \TrS b c^2 x^2 \ y \TrS b^2 x^3 y \TrS a b^3 x^3 y \TrS b^2
c x^3 y \TrS a b^3 c x^3 y \TrS b^2 c^2 x^3 y \TrS a^2 \ y^2 \TrS
a^3 b y^2 \TrS a^2 c y^2 \TrS a x y^2 \TrS a^2 b x y^2 \TrS a^3
b^2 x y^2 \TrS a c x \ y^2 \TrS a^2 b c x y^2 \TrS a c^2 x y^2
\TrS x^2 y^2 \TrS a b x^2 y^2 \TrS a^2 b^2 x^2 y^2 \TrS c \ x^2
y^2 \TrS a b c x^2 y^2 \TrS a^2 b^2 c x^2 y^2 \TrS c^2 x^2 y^2
\TrS a b c^2 x^2 y^2 \TrS b \ x^3 y^2 \TrS a b^2 x^3 y^2 \TrS b c
x^3 y^2 \TrS a b^2 c x^3 y^2 \TrS a b^2 c^2 x^3 \ y^2 \TrS a^3 y^3
\TrS a^3 c y^3 \TrS a^2 x y^3 \TrS a^3 b x y^3 \TrS a^2 c x y^3
\TrS a^3 b c x \ y^3 \TrS a^2 c^2 x y^3 \TrS a x^2 y^3 \TrS a^2 b
x^2 y^3 \TrS a c x^2 y^3 \TrS a^2 b c x^2 \ y^3 \TrS a^2 b c^2 x^2
y^3 \TrS a b z \TrS a b c z \TrS b x z \TrS a b^2 x z \TrS b c x z
\TrS a b^2 c x \ z \TrS b c^2 x z \TrS b^2 x^2 z \TrS a b^3 x^2 z
\TrS b^2 c x^2 z \TrS b^2 c^2 x^2 z \TrS b^3 x^3 \ z \TrS b^3 c
x^3 z \TrS a y z \TrS a^2 b y z \TrS a c y z \TrS a^2 b c y z \TrS
a c^2 y z \TrS a b x y \ z \TrS a^2 b^2 x y z \TrS c x y z \TrS a
b c x y z \TrS a^2 b^2 c x y z \TrS c^2 x y z \TrS a b c^2 \ x y z
\TrS c^3 x y z \TrS b x^2 y z \TrS a b^2 x^2 y z \TrS a^2 b^3 x^2
y z \TrS b c x^2 y z \TrS a \ b^2 c x^2 y z \TrS b c^2 x^2 y z
\TrS a b^2 c^2 x^2 y z \TrS b c^3 x^2 y z \TrS b^2 x^3 y \ z \TrS
a b^3 x^3 y z \TrS b^2 c x^3 y z \TrS a b^3 c x^3 y z \TrS b^2 c^2
x^3 y z \TrS a^2 y^2 \ z \TrS a^3 b y^2 z \TrS a^2 c y^2 z \TrS
a^2 c^2 y^2 z \TrS a x y^2 z \TrS a^2 b x y^2 z \TrS a^3 \ b^2 x
y^2 z \TrS a c x y^2 z \TrS a^2 b c x y^2 z \TrS a c^2 x y^2 z
\TrS a^2 b c^2 x y^2 \ z \TrS a c^3 x y^2 z \TrS a b x^2 y^2 z
\TrS a^2 b^2 x^2 y^2 z \TrS c x^2 y^2 z \TrS a b c x^2 \ y^2 z
\TrS a^2 b^2 c x^2 y^2 z \TrS c^2 x^2 y^2 z \TrS a b c^2 x^2 y^2 z
\TrS a b c^3 x^2 \ y^2 z \TrS a^3 y^3 z \TrS a^3 c y^3 z \TrS a^2
x y^3 z \TrS a^3 b x y^3 z \TrS a^2 c x y^3 \ z \TrS a^3 b c x y^3
z \TrS a^2 c^2 x y^3 z \TrS a b z^2 \TrS a^2 b^2 z^2 \TrS a b c
z^2 \TrS a b^2 \ x z^2 \TrS b c x z^2 \TrS a b^2 c x z^2 \TrS
{ \bf b c^2 x z^2}
 \TrS b^2 x^2 z^2 \TrS a b^3 x^2 \ z^2 \TrS b^2 c x^2 z^2 \TrS a b^3 c x^2 z^2 \TrS b^2 c^2
x^2 z^2 \TrS b^3 x^3 z^2 \TrS b^3 c \ x^3 z^2 \TrS a^2 b y z^2
\TrS a c y z^2 \TrS a^2 b c y z^2 \TrS
{\bf  a c^2 y z^2}
 \TrS a b x y \ z^2 \TrS a^2 b^2 x y z^2 \TrS a b c x y z^2 \TrS a^2 b^2 c x y z^2 \TrS a b
c^2 x y \ z^2 \TrS c^3 x y z^2 \TrS a b^2 x^2 y z^2 \TrS a^2 b^3
x^2 y z^2 \TrS b c x^2 y z^2 \TrS a b^2 \ c x^2 y z^2 \TrS b c^2
x^2 y z^2 \TrS a b^2 c^2 x^2 y z^2 \TrS b c^3 x^2 y z^2 \TrS a^2 \
y^2 z^2 \TrS a^3 b y^2 z^2 \TrS a^2 c y^2 z^2 \TrS a^3 b c y^2 z^2
\TrS a^2 c^2 y^2 \ z^2 \TrS a^2 b x y^2 z^2 \TrS a^3 b^2 x y^2 z^2
\TrS a c x y^2 z^2 \TrS a^2 b c x y^2 \ z^2 \TrS a c^2 x y^2 z^2
\TrS a^2 b c^2 x y^2 z^2 \TrS a c^3 x y^2 z^2 \TrS a^3 y^3 \ z^2
\TrS a^3 c y^3 z^2 \TrS a^2 b^2 z^3 \TrS a b c z^3 \TrS a b c^2
z^3 \TrS a b^2 x z^3 \TrS a^2 \ b^3 x z^3 \TrS a b^2 c x z^3 \TrS
a b^2 c^2 x z^3 \TrS b c^3 x z^3 \TrS a b^3 x^2 z^3 \TrS b^2 \ c
x^2 z^3 \TrS a b^3 c x^2 z^3 \TrS b^2 c^2 x^2 z^3 \TrS a^2 b y z^3
\TrS a^3 b^2 y \ z^3 \TrS a^2 b c y z^3 \TrS a^2 b c^2 y z^3 \TrS
a c^3 y z^3 \TrS a^2 b^2 x y z^3 \TrS a b c x \ y z^3 \TrS a^2 b^2
c x y z^3 \TrS a b c^2 x y z^3 \TrS c^3 x y z^3 \TrS a b c^3 x y \
z^3 \TrS a^3 b y^2 z^3 \TrS a^2 c y^2 z^3 \TrS a^3 b c y^2 z^3
\TrS a^2 c^2 y^2 z^3 \TrS a^2 \ b^2 z^4 \TrS a^2 b^2 c z^4 \TrS
a^2 b^3 x z^4 \TrS a b^2 c x z^4 \TrS a b^2 c^2 x z^4 \TrS b \ c^3
x z^4 \TrS a^3 b^2 y z^4 \TrS a^2 b c y z^4 \TrS a^2 b c^2 y z^4
\TrS a c^3 y z^4 \TrS a^2 \ b^2 c z^5 \TrS a b c^3 z^5$ \pSkip

\noindent  $\al_{11} =a^3 b^2 x^2 y^3 \TrS a c^2 x^2 y^3 \TrS a^3
b^2 y z^2$
 \pSkip

\noindent $ \bt_{11}= a^2 b^3 x^3 y^2 \TrS b c^2 x^3 y^2 \TrS a^2
b^3 x z^2 $ \pSkip

By Equation  \Ref{eq:a1} there are two dominate monomials for
which: \pSkip
 \noindent  $\al_{11} =a^3 b^2 x^2 y^3 \TrS a c^2 x^2 y^3 \TrS a^3 b^2 y
z^2 \leq { \bf a c^2 y z^2}$ \pSkip
\noindent $ \bt_{11}= a^2 b^3 x^3 y^2 \TrS b c^2 x^3 y^2 \TrS a^2
b^3 x z^2 \leq { \bf b c^2 x z^2}$ .\pSkip
Thus, $\al_{11}$ and $\bt_{11}$  are inessential for $f_{11}$ and
$g_{11}$ respectively, namely $f_{11} \eqR g_{11}$. \pSkip

 \pSkip


 \item[$\underline{i=1, j =2}$]  $ $ \pSkip

\noindent $h_{12} =  \overline{a}^{(\bt.3.)} \TrS a^2 b \TrS a b x
\TrS a^2 b^2 x \TrS
\overline{c x}^{(\bt.4.)}
 \TrS a b c x \TrS a b^2 x^2 \TrS b c x^2 \TrS a b^2 c x^2 \TrS b c^2 x^2 \TrS b^2 c x^3 \TrS b^2 c^2
x^3 \TrS a^2 y \TrS a^3 b y \TrS a^2 c y \TrS a x y \TrS a^2 b x y
\TrS
\overline{\underline{a^3 b^2 x y}}_{(\al.2.)}^{(\bt.5.)}
 \TrS a c x y \TrS a^2 b c x y \TrS a c^2 x y \TrS a b x^2 y \TrS a^2 b^2 x^2 y \TrS
\overline{c x^2 y}^{(\bt.2.)}
 \TrS a b c x^2 y \TrS a^2 b^2 c x^2 y \TrS c^2 x^2 y \TrS a b c^2 x^2 y \TrS
\underline{c^3 x^2 y}_{(\al.3.)}
 \TrS a b^2 x^3 y \TrS b c x^3 y \TrS a b^2 c x^3 y \TrS b c^2 x^3 y \TrS a b^2 c^2 x^3 y \TrS b
c^3 x^3 y \TrS a^3 y^2 \TrS a^4 b y^2 \TrS a^3 c y^2 \TrS a^2 x
y^2 \TrS a^3 b x y^2 \TrS a^2 c x y^2 \TrS a^3 b c x y^2 \TrS a^2
c^2 x y^2 \TrS a x^2 y^2 \TrS a^2 b x^2 y^2 \TrS a^3 b^2 x^2 y^2
\TrS a c x^2 y^2 \TrS a^2 b c x^2 y^2 \TrS a c^2 x^2 y^2 \TrS a^2
b c^2 x^2 y^2 \TrS a c^3 x^2 y^2 \TrS
\overline{\underline{a b x^3 y^2}}_{(\al.1.)}^{(\bt.1.)}
 \TrS a^2 b^2 x^3 y^2 \TrS a^2 b^2 c x^3 y^2 \TrS
\overline{\underline{a b c^3 x^3 y^2}}_{(\al.1.)}^{(\bt.1.)}
 \TrS a^4 y^3 \TrS a^4 c y^3 \TrS a^3 x y^3 \TrS a^4 b x y^3 \TrS a^3 c x y^3 \TrS a^3 c^2 x
y^3 \TrS a^2 x^2 y^3 \TrS a^3 b x^2 y^3 \TrS a^2 c x^2 y^3 \TrS
a^3 b c x^2 y^3 \TrS a^2 c^2 x^2 y^3 \TrS a^2 c^3 x^2 y^3 \TrS
\underline{a^2 b z}_{(\al.4.)}
 \TrS a^2 b c z \TrS a b x z \TrS a^2 b^2 x z \TrS
\underline{c x z}_{(\al.5.)}
 \TrS a b c x z \TrS a b c^2 x z \TrS a b^2 x^2 z \TrS b c x^2 z \TrS a b^2 c x^2 z \TrS b c^2
x^2 z \TrS b c^3 x^2 z \TrS b^2 c x^3 z \TrS b^2 c^2 x^3 z \TrS
a^2 y z \TrS a^3 b y z \TrS a^2 c y z \TrS a^3 b c y z \TrS a^2
c^2 y z \TrS a^2 b x y z \TrS a c x y z \TrS a^2 b c x y z \TrS a
c^2 x y z \TrS a^2 b c^2 x y z \TrS a c^3 x y z \TrS a b x^2 y z
\TrS a^2 b^2 x^2 y z \TrS c x^2 y z \TrS a b c x^2 y z \TrS a^2
b^2 c x^2 y z \TrS a b c^2 x^2 y z \TrS a b c^3 x^2 y z \TrS c^4
x^2 y z \TrS a b^2 x^3 y z \TrS b c x^3 y z \TrS a b^2 c x^3 y z
\TrS b c^2 x^3 y z \TrS a b^2 c^2 x^3 y z \TrS b c^3 x^3 y z \TrS
a^3 y^2 z \TrS a^4 b y^2 z \TrS a^3 c y^2 z \TrS a^3 c^2 y^2 z
\TrS a^2 x y^2 z \TrS a^3 b x y^2 z \TrS a^2 c x y^2 z \TrS a^3 b
c x y^2 z \TrS a^2 c^2 x y^2 z \TrS a^2 c^3 x y^2 z \TrS a^2 b x^2
y^2 z \TrS a^3 b^2 x^2 y^2 z \TrS a c x^2 y^2 z \TrS a^2 b c x^2
y^2 z \TrS a c^2 x^2 y^2 z \TrS a^2 b c^2 x^2 y^2 z \TrS a c^3 x^2
y^2 z \TrS a c^4 x^2 y^2 z \TrS a^4 y^3 z \TrS a^4 c y^3 z \TrS
a^3 x y^3 z \TrS a^4 b x y^3 z \TrS a^3 c x y^3 z \TrS a^3 c^2 x
y^3 z \TrS
\overline{a^2 b z^2}^{(\bt.6.)}
 \TrS \underline{a c z^2}_{(\al.6.)}
 \TrS a^2 b^2 x z^2 \TrS a b c x z^2 \TrS a^2 b^2 c x z^2 \TrS a b c^2 x z^2 \TrS b c x^2
z^2 \TrS a b^2 c x^2 z^2 \TrS b c^2 x^2 z^2 \TrS a b^2 c^2 x^2 z^2
\TrS b c^3 x^2 z^2 \TrS b^2 c x^3 z^2 \TrS b^2 c^2 x^3 z^2 \TrS
a^3 b y z^2 \TrS a^2 c y z^2 \TrS a^3 b c y z^2 \TrS a^2 c^2 y z^2
\TrS a^2 b x y z^2 \TrS a c x y z^2 \TrS a^2 b c x y z^2 \TrS a
c^2 x y z^2 \TrS a^2 b c^2 x y z^2 \TrS a c^3 x y z^2 \TrS a^2 b^2
x^2 y z^2 \TrS a b c x^2 y z^2 \TrS a^2 b^2 c x^2 y z^2 \TrS c^2
x^2 y z^2 \TrS a b c^2 x^2 y z^2 \TrS
\overline{\underline{c^3 x^2 y z^2}}_{(\al.3.)}^{(\bt.2.)}
 \TrS a b c^3 x^2 y z^2 \TrS c^4 x^2 y z^2 \TrS a^3 y^2 z^2 \TrS a^4 b y^2 z^2 \TrS a^3 c
y^2 z^2 \TrS a^3 c^2 y^2 z^2 \TrS a^3 b x y^2 z^2 \TrS a^2 c x y^2
z^2 \TrS a^3 b c x y^2 z^2 \TrS a^2 c^2 x y^2 z^2 \TrS a^2 c^3 x
y^2 z^2 \TrS a^4 y^3 z^2 \TrS a^4 c y^3 z^2 \TrS
\underline{a^2 b c^2 z^3}_{(\al.4.)}
 \TrS  \overline{a c^3 z^3}^{(\bt.3.)}
 \TrS a b c x z^3 \TrS a^2 b^2 c x z^3 \TrS a b c^2 x z^3 \TrS a b c^3 x z^3 \TrS
\overline{c^4 x z^3}^{(\bt.4.)}
 \TrS a b^2 c x^2 z^3 \TrS b c^2 x^2 z^3 \TrS a b^2 c^2 x^2 z^3 \TrS b c^3 x^2 z^3 \TrS a^3
b y z^3 \TrS a^2 c y z^3 \TrS a^3 b c y z^3 \TrS a^2 c^2 y z^3
\TrS a^2 c^3 y z^3 \TrS
\overline{\underline{a^3 b^2 x y z^3}}_{(\al.2.)}^{(\bt.5.)}
 \TrS a^2 b c x y z^3 \TrS a c^2 x y z^3 \TrS a^2 b c^2 x y z^3 \TrS a c^3 x y z^3 \TrS a
c^4 x y z^3 \TrS a^4 b y^2 z^3 \TrS a^3 c y^2 z^3 \TrS a^3 c^2 y^2
z^3 \TrS a^2 b c z^4 \TrS
\overline{a^2 b c^2 z^4}^{(\bt.6.)}
 \TrS a^2 b^2 c x z^4 \TrS a b c^2 x z^4 \TrS a b c^3 x z^4 \TrS
\underline{c^4 x z^4}_{(\al.5.)}
 \TrS a^3 b c y z^4 \TrS a^2 c^2 y z^4 \TrS a^2 c^3 y z^4 \TrS a^2 b c^2 z^5 \TrS
\underline{a c^4 z^5}_{(\al.6.)}
$ \pSkip

\noindent $\al_{12} = c^2 x^3 y^2 \TrS a b c^2 x^3 y^2 \TrS a^3
b^2 x y z \TrS c^3 x^2 y z \TrS a^2 b c z^2 \TrS { \bf c^2 x z^2 }
 \TrS a^3 b^2 z^3 \TrS { \bf a c^2 z^3}$ \pSkip

\noindent $ \bt_{12} = a b c x^3 y^2 \TrS c^3 x^3 y^2 \TrS c^2 x^2
y z \TrS a^3 b^2 z^2 \TrS { \bf a c^2 z^2 \TrS   c^3 x \ z^2} \TrS
a^3 b^2 x y z^2 \TrS a^2 b c z^3$ \pSkip

Inside $\al_{12}$, using Equation  \Ref{eq:a1}, we have $ c^2 x^3
y^2 \leq {\bf c^2 x z^2}$ and $a^3 b^2 z^3 \leq {\bf a c^2 z^3}$,
so they are inessential. By the following diagram we show that the
other monomials of $\al_{12}$ are inessential with respect to
$h_{12}$.   \pSkip

\noindent
$$  \begin{array}{l} \\ \\  \al_{12} = \end{array}
\begin{array}{cccccccccccc}
         \al.1. & & \al.2. & & \al.3. & & \al.4. & & \al.5. & &  \al.6.
         \\  \hline &&&&& &&&&&\\
         a b x^3 y^2 & & a^3 b^2 x y & & c^3 x^2 y  & &  a^2 b z & & c x z & & a c z^2 \\
         | & &  |& & |& & |& & |& & | \\
  a b c^2 x^3 y^2 & \TrS  & a^3 b^2 x y z &  \TrS  & c^3 x^2 y
z& \TrS & a^2 b c z^2 &  \TrS  & c^2 x z^2 & \TrS & a c^2 z^3  \\
   | & &  |& & |& & |& & |& & | \\
  a b c^3 x^3 y^2 & &  a^3 b^2 x y z^3 & & c^3 x^2 y z^2 & & a^2 b c^2 z^3 & &
  c^4 x z^4 & & a c^4 z^5 \\
\end{array}
$$
The upper and the lower rows specify the monomials of $h_{12}$
that make respectively each term of $\al_{12}$  to be inessential.
The monomials of $h_{12}$ are correspondingly marked above.

Inside $\bt_{12}$, again using Equation  \Ref{eq:a1}, we have $
c^3 x^3 y^2 \leq {\bf c^3 x \ z^2}$ and $a^3 b^2 z^2 \leq {\bf a
c^2 z^2 }$, so these monomials are inessential. The below diagram
shows that all the other monomials of $\bt_{12}$ are inessential
with respect to $h_{12}$.

\noindent
$$  \begin{array}{l} \\ \\  \bt _{12} = \end{array}
\begin{array}{cccccccccccc}
         \bt.1. & & \bt.2. && \bt.3. & & \bt.4. &&  \bt.5. && \bt.6.
         \\  \hline &&&&& &&&&&\\
        a b x^3 y^2 & & c x^2 y & & a && cx & &a^3 b^2 x y&& a^2 b z^2\\
          | & &  |& & |& & |& & |& & | \\
   a b c x^3 y^2 &  \TrS  & c^2 x^2 y z &  \TrS  & a c^2 z^2 &  \TrS  & c^3 x \
z^2 &  \TrS  & a^3 b^2 x y z^2 &  \TrS  & a^2 b c z^3  \\
   | & &  |& & |& & |& & |& & | \\
a b c^3 x^3 y^2 & & c^3 x^2 y z^2 & & a c^3 z^3&& c^4 x z^3 && a^3 b^2 x y z^3 && a^2 b c^2 z^4\\
\end{array}
$$
%
Accordingly,  $f_{12} \eqR g_{12}$.

\pSkip
 \item[$\underline{i=2, j =1}$] This case is parallel to the case
 when $i=1$ and $j = 2$.
  \pSkip
\noindent $h_{21} = \underline{b}_{(\al.3.)}  \TrS a b^2 \TrS b^2
x \TrS a b^3 x \TrS b^2 c x \TrS b^3 x^2 \TrS a b^4 x^2 \TrS b^3 c
x^2 \TrS b^4 x^3 \TrS b^4 c x^3 \TrS a b y \TrS a^2 b^2 y \TrS
\underline{c y}_{(\al.4.)}
\TrS a b c y \TrS b x y \TrS a b^2 x y \TrS
 \overline{\underline{a^2 b^3 x y}}_{(\al.5.)}^{(\bt.2.)} \TrS
 b c
x y \TrS a b^2 c x y \TrS b c^2 x y \TrS b^2 x^2 y \TrS a b^3 x^2
y \TrS b^2 c x^2 y \TrS a b^3 c x^2 y \TrS b^2 c^2 x^2 y \TrS b^3
x^3 y \TrS a b^4 x^3 y \TrS b^3 c x^3 y \TrS b^3 c^2 x^3 y \TrS
a^2 b y^2 \TrS a c y^2 \TrS a^2 b c y^2 \TrS a c^2 y^2 \TrS a b x
y^2 \TrS a^2 b^2 x y^2 \TrS
\underline{c x y^2}_{(\al.2.)}
\TrS a b c x y^2 \TrS a^2 b^2 c x y^2 \TrS c^2 x y^2 \TrS a b c^2
x y^2 \TrS
\overline{c^3 x y^2}^{(\bt.3.)}
 \TrS b x^2 y^2 \TrS a b^2 x^2 y^2 \TrS a^2 b^3 x^2 y^2
\TrS b c x^2 y^2 \TrS a b^2 c x^2 y^2 \TrS b c^2 x^2 y^2 \TrS a
b^2 c^2 x^2 y^2 \TrS b c^3 x^2 y^2 \TrS b^2 x^3 y^2 \TrS a b^3 x^3
y^2 \TrS b^2 c x^3 y^2 \TrS a b^3 c x^3 y^2 \TrS b^2 c^2 x^3 y^2
\TrS b^2 c^3 x^3 y^2 \TrS a^2 c y^3 \TrS a^2 c^2 y^3 \TrS a^2 b x
y^3 \TrS a c x y^3 \TrS a^2 b c x y^3 \TrS a c^2 x y^3 \TrS a^2 b
c^2 x y^3 \TrS a c^3 x y^3 \TrS
\overline{\underline{ a b x^2 y^3}}_{(\al.1.)}^{(\bt.1.)} \TrS
a^2 b^2 x^2 y^3 \TrS a^2 b^2 c x^2 y^3 \TrS
 \overline{\underline{a b c^3 x^2 y^3}}_{(\al.1.)}^{(\bt.1.)}
 \TrS  \overline{a b^2 z}^{(\bt.4.)}
 \TrS a b^2 c z \TrS b^2 x z \TrS a b^3 x z \TrS b^2 c x z
\TrS a b^3 c x z \TrS b^2 c^2 x z \TrS b^3 x^2 z \TrS a b^4 x^2 z
\TrS b^3 c x^2 z \TrS b^3 c^2 x^2 z \TrS b^4 x^3 z \TrS b^4 c x^3
z \TrS a b y z \TrS a^2 b^2 y z \TrS
\overline{c y z}^{(\bt.5.)}
\TrS a b c y z \TrS a b c^2 y z \TrS a b^2 x y z \TrS b c x y z
\TrS a b^2 c x y z \TrS b c^2 x y z \TrS a b^2 c^2 x y z \TrS b
c^3 x y z \TrS b^2 x^2 y z \TrS a b^3 x^2 y z \TrS b^2 c x^2 y z
\TrS a b^3 c x^2 y z \TrS b^2 c^2 x^2 y z \TrS b^2 c^3 x^2 y z
\TrS b^3 x^3 y z \TrS a b^4 x^3 y z \TrS b^3 c x^3 y z \TrS b^3
c^2 x^3 y z \TrS a^2 b y^2 z \TrS a c y^2 z \TrS a^2 b c y^2 z
\TrS a c^2 y^2 z \TrS a c^3 y^2 z \TrS a b x y^2 z \TrS a^2 b^2 x
y^2 z \TrS c x y^2 z \TrS a b c x y^2 z \TrS a^2 b^2 c x y^2 z
\TrS a b c^2 x y^2 z \TrS a b c^3 x y^2 z \TrS c^4 x y^2 z \TrS a
b^2 x^2 y^2 z \TrS a^2 b^3 x^2 y^2 z \TrS b c x^2 y^2 z \TrS a b^2
c x^2 y^2 z \TrS b c^2 x^2 y^2 z \TrS a b^2 c^2 x^2 y^2 z \TrS b
c^3 x^2 y^2 z \TrS b c^4 x^2 y^2 z \TrS a^2 c y^3 z \TrS a^2 c^2
y^3 z \TrS a^2 b x y^3 z \TrS a c x y^3 z \TrS a^2 b c x y^3 z
\TrS a c^2 x y^3 z \TrS a^2 b c^2 x y^3 z \TrS a c^3 x y^3 z \TrS
\underline{a b^2 z^2}_{(\al.6)}
\TrS\overline{ b c z^2}^{(\bt.6.)}
\TrS a b^3 x z^2 \TrS b^2 c x z^2 \TrS a b^3 c x z^2 \TrS b^2 c^2
x z^2 \TrS b^3 x^2 z^2 \TrS a b^4 x^2 z^2 \TrS b^3 c x^2 z^2 \TrS
b^3 c^2 x^2 z^2 \TrS b^4 x^3 z^2 \TrS b^4 c x^3 z^2 \TrS a^2 b^2 y
z^2 \TrS a b c y z^2 \TrS a^2 b^2 c y z^2 \TrS a b c^2 y z^2 \TrS
a b^2 x y z^2 \TrS b c x y z^2 \TrS a b^2 c x y z^2 \TrS b c^2 x y
z^2 \TrS a b^2 c^2 x y z^2 \TrS b c^3 x y z^2 \TrS a b^3 x^2 y z^2
\TrS b^2 c x^2 y z^2 \TrS a b^3 c x^2 y z^2 \TrS b^2 c^2 x^2 y z^2
\TrS b^2 c^3 x^2 y z^2 \TrS a c y^2 z^2 \TrS a^2 b c y^2 z^2 \TrS
a c^2 y^2 z^2 \TrS a^2 b c^2 y^2 z^2 \TrS a c^3 y^2 z^2 \TrS a^2
b^2 x y^2 z^2 \TrS a b c x y^2 z^2 \TrS a^2 b^2 c x y^2 z^2 \TrS
c^2 x y^2 z^2 \TrS a b c^2 x y^2 z^2 \TrS
\overline{\underline{c^3 x y^2 z^2 }}_{(\al.2.)}^{(\bt.3.)}
\TrS a b c^3 x y^2 z^2 \TrS c^4 x y^2 z^2 \TrS a^2 c y^3 z^2 \TrS
a^2 c^2 y^3 z^2 \TrS
\overline{a b^2 c^2 z^3}^{(\bt.4.)} \TrS
 \underline{b c^3 z^3}_{(\al.3.)}
 \TrS a b^3 x z^3
\TrS b^2 c x z^3 \TrS a b^3 c x z^3 \TrS b^2 c^2 x z^3 \TrS b^2
c^3 x z^3 \TrS a b^4 x^2 z^3 \TrS b^3 c x^2 z^3 \TrS b^3 c^2 x^2
z^3 \TrS a b c y z^3 \TrS a^2 b^2 c y z^3 \TrS a b c^2 y z^3 \TrS
a b c^3 y z^3 \TrS
\underline{ c^4 y z^3}_{(\al.4.)}
 \TrS \overline{\underline{a^2 b^3 x y z^3}}_{(\al.5.)}^{(\bt.2.)}
 \TrS a b^2 c x y z^3 \TrS b c^2 x y z^3 \TrS a b^2
c^2 x y z^3 \TrS b c^3 x y z^3 \TrS b c^4 x y z^3 \TrS a^2 b c y^2
z^3 \TrS a c^2 y^2 z^3 \TrS a^2 b c^2 y^2 z^3 \TrS a c^3 y^2 z^3
\TrS a b^2 c z^4 \TrS
\underline{a b^2 c^2 z^4}_{(\al.6.)}
\TrS a b^3 c x z^4 \TrS b^2 c^2 x z^4 \TrS b^2 c^3 x z^4 \TrS a^2
b^2 c y z^4 \TrS a b c^2 y z^4 \TrS a b c^3 y z^4 \TrS
\overline{c^4 y z^4}^{(\bt.5.)}
\TrS a b^2 c^2 z^5 \TrS
\overline{b c^4 z^5}^{(\bt.6.)}$

\pSkip

\noindent  $\al_{21} =a b c x^2 y^3 \TrS c^3 x^2 y^3 \TrS c^2 x
y^2 z \TrS a^2 b^3 z^2 \TrS {\bf b c^2 z^2} \TrS {\bf c^3 y z^2}
\TrS a^2 b^3 x y z^2 \TrS a b^2 c z^3$ \pSkip

\noindent $ \bt_{21} = c^2 x^2 y^3 \TrS a b c^2 x^2 y^3 \TrS a^2
b^3 x y z \TrS c^3 x y^2 z \TrS a b^2 c z^2 \TrS {\bf c^2 y z^2}
\TrS a^2 b^3 z^3 \TrS {\bf b c^2 z^3}$ \pSkip

Inside $\al_{21}$, again using Equation  \Ref{eq:a1}, we have $
c^3 x^2 y^3 \leq {\bf c^3 y \ z^2}$ and $a^2 b^3 z^2 \leq {\bf b
c^2 z^2 }$, so these monomials are inessential. The below diagram
shows that all the other monomials of $\al_{21}$ are inessential
with respect to $h_{21}$.

\noindent
$$  \begin{array}{l} \\ \\  \al _{21} = \end{array}
\begin{array}{cccccccccccc}
         \al.1. & & \al.2. & & \al.3. & & \al.4. & & \al.5. & &  \al.6.
         \\  \hline &&&&& &&&&&\\
        a b x^2 y^3 & & c x y^2 & & b && cy & &a^2 b^3 x y&& a b^2 z^2\\
          | & &  |& & |& & |& & |& & | \\
   a b c x^2 y^3 &  \TrS  & c^2 x y^2 z &  \TrS  & b c^2 z^2 &  \TrS  & c^3 y \
z^2 &  \TrS  & a^2 b^3 x y z^2 &  \TrS  & a b^2 c z^3  \\
   | & &  |& & |& & |& & |& & | \\
a b c^3 x^2 y^3 & & c^3 x y^2 z^2 & & b c^3 z^3&& c^4 y z^3 && a^2 b^3 x y z^3 && a b^2 c^2 z^4\\
\end{array}
$$

Inside $\bt_{21}$, using Equation  \Ref{eq:a1}, we have $ c^2 x^2
y^3 \leq {\bf c^2 y z^2}$ and $a^2 b^3 z^3 \leq {\bf b c^2 z^3}$,
so they are inessential. By the following diagram we show that the
other monomials of $\bt_{21}$ are inessential with respect to
$h_{21}$.   \pSkip

\noindent
$$  \begin{array}{l} \\ \\  \bt_{12} = \end{array}
\begin{array}{cccccccccccc}
\bt.1. & & \bt.2. && \bt.3. & & \bt.4. &&  \bt.5. && \bt.6.
         \\  \hline &&&&& &&&&&\\
         a b x^2 y^3 & & a^2 b^3 x y & & c^3 x y^2  & &  a b^2 z & & c y z & & b c z^2 \\
         | & &  |& & |& & |& & |& & | \\
  a b c^2 x^2 y^3 & \TrS  & a^2 b^3 x y z &  \TrS  & c^3 x y^2
z& \TrS & a b^2 c z^2 &  \TrS  & c^2 y z^2 & \TrS & b c^2 z^3  \\
   | & &  |& & |& & |& & |& & | \\
  a b c^3 x^2 y^3 & &  a^2 b^3 x y z^3 & & c^3 x y^2 z^2 & & a b^2 c^2 z^3 & &
  c^4 y z^4 & & b c^4 z^5 \\
\end{array}
$$

Accordingly, $f_{21} \eqR g_{21}$.

\pSkip

  \item[$\underline{i=2, j =2}$] This case is parallel to the case when $i=j=1$.

  \pSkip

\noindent $h_{22} = a b \TrS a^2 b^2 \TrS a b^2 x \TrS a^2 b^3 x
\TrS b c x \TrS a b^2 c x \TrS a b^3 x^2 \TrS b^2 c x^2 \TrS a b^3
c x^2 \TrS b^2 c^2 x^2 \TrS b^3 c x^3 \TrS b^3 c^2 x^3 \TrS a^2 b
y \TrS a^3 b^2 y \TrS a c y \TrS a^2 b c y \TrS a b x y \TrS a^2
b^2 x y \TrS a b c x y \TrS a^2 b^2 c x y \TrS c^2 x y \TrS a b
c^2 x y \TrS a b^2 x^2 y \TrS a^2 b^3 x^2 y \TrS b c x^2 y \TrS a
b^2 c x^2 y \TrS b c^2 x^2 y \TrS a b^2 c^2 x^2 y \TrS b c^3 x^2 y
\TrS a b^3 x^3 y \TrS b^2 c x^3 y \TrS a b^3 c x^3 y \TrS b^2 c^2
x^3 y \TrS b^2 c^3 x^3 y \TrS a^3 b y^2 \TrS a^2 c y^2 \TrS a^3 b
c y^2 \TrS a^2 c^2 y^2 \TrS a^2 b x y^2 \TrS a^3 b^2 x y^2 \TrS a
c x y^2 \TrS a^2 b c x y^2 \TrS a c^2 x y^2 \TrS a^2 b c^2 x y^2
\TrS a c^3 x y^2 \TrS a b x^2 y^2 \TrS a^2 b^2 x^2 y^2 \TrS a b c
x^2 y^2 \TrS a^2 b^2 c x^2 y^2 \TrS a b c^2 x^2 y^2 \TrS c^3 x^2
y^2 \TrS a b c^3 x^2 y^2 \TrS c^4 x^2 y^2 \TrS a b^2 x^3 y^2 \TrS
a b^2 c x^3 y^2 \TrS a b^2 c^2 x^3 y^2 \TrS b c^3 x^3 y^2 \TrS b
c^4 x^3 y^2 \TrS a^3 c y^3 \TrS a^3 c^2 y^3 \TrS a^3 b x y^3 \TrS
a^2 c x y^3 \TrS a^3 b c x y^3 \TrS a^2 c^2 x y^3 \TrS a^2 c^3 x
y^3 \TrS a^2 b x^2 y^3 \TrS a^2 b c x^2 y^3 \TrS a^2 b c^2 x^2 y^3
\TrS a c^3 x^2 y^3 \TrS a c^4 x^2 y^3 \TrS a^2 b^2 z \TrS a^2 b^2
c z \TrS a b^2 x z \TrS a^2 b^3 x z \TrS b c x z \TrS a b^2 c x z
\TrS a b^2 c^2 x z \TrS a b^3 x^2 z \TrS b^2 c x^2 z \TrS a b^3 c
x^2 z \TrS b^2 c^2 x^2 z \TrS b^2 c^3 x^2 z \TrS b^3 c x^3 z \TrS
b^3 c^2 x^3 z \TrS a^2 b y z \TrS a^3 b^2 y z \TrS a c y z \TrS
a^2 b c y z \TrS a^2 b c^2 y z \TrS a^2 b^2 x y z \TrS a b c x y z
\TrS a^2 b^2 c x y z \TrS c^2 x y z \TrS a b c^2 x y z \TrS a b
c^3 x y z \TrS a b^2 x^2 y z \TrS a^2 b^3 x^2 y z \TrS b c x^2 y z
\TrS a b^2 c x^2 y z \TrS b c^2 x^2 y z \TrS a b^2 c^2 x^2 y z
\TrS b c^3 x^2 y z \TrS b c^4 x^2 y z \TrS a b^3 x^3 y z \TrS b^2
c x^3 y z \TrS a b^3 c x^3 y z \TrS b^2 c^2 x^3 y z \TrS b^2 c^3
x^3 y z \TrS a^3 b y^2 z \TrS a^2 c y^2 z \TrS a^3 b c y^2 z \TrS
a^2 c^2 y^2 z \TrS a^2 c^3 y^2 z \TrS a^2 b x y^2 z \TrS a^3 b^2 x
y^2 z \TrS a c x y^2 z \TrS a^2 b c x y^2 z \TrS a c^2 x y^2 z
\TrS a^2 b c^2 x y^2 z \TrS a c^3 x y^2 z \TrS a c^4 x y^2 z \TrS
a^2 b^2 x^2 y^2 z \TrS a b c x^2 y^2 z \TrS a^2 b^2 c x^2 y^2 z
\TrS a b c^2 x^2 y^2 z \TrS c^3 x^2 y^2 z \TrS a b c^3 x^2 y^2 z
\TrS c^4 x^2 y^2 z \TrS c^5 x^2 y^2 z \TrS a^3 c y^3 z \TrS a^3
c^2 y^3 z \TrS a^3 b x y^3 z \TrS a^2 c x y^3 z \TrS a^3 b c x y^3
z \TrS a^2 c^2 x y^3 z \TrS a^2 c^3 x y^3 z \TrS a b c z^2 \TrS
a^2 b^2 c z^2 \TrS a b c^2 z^2 \TrS a b^2 c x z^2 \TrS
{ \bf b c^2 x z^2 }
 \TrS a b^2 c^2 x z^2 \TrS b c^3 x z^2 \TrS b^2 c x^2 z^2 \TrS a b^3 c x^2 z^2 \TrS b^2 c^2
x^2 z^2 \TrS b^2 c^3 x^2 z^2 \TrS b^3 c x^3 z^2 \TrS b^3 c^2 x^3
z^2 \TrS a^2 b c y z^2 \TrS
{\bf  a c^2 y z^2}
 \TrS a^2 b c^2 y z^2 \TrS a c^3 y z^2 \TrS a^2 b^2 x y z^2 \TrS a b c x y z^2 \TrS a^2 b^2
c x y z^2 \TrS c^2 x y z^2 \TrS a b c^2 x y z^2 \TrS c^3 x y z^2
\TrS a b c^3 x y z^2 \TrS c^4 x y z^2 \TrS a^2 b^3 x^2 y z^2 \TrS
a b^2 c x^2 y z^2 \TrS b c^2 x^2 y z^2 \TrS a b^2 c^2 x^2 y z^2
\TrS b c^3 x^2 y z^2 \TrS b c^4 x^2 y z^2 \TrS a^2 c y^2 z^2 \TrS
a^3 b c y^2 z^2 \TrS a^2 c^2 y^2 z^2 \TrS a^2 c^3 y^2 z^2 \TrS a^3
b^2 x y^2 z^2 \TrS a^2 b c x y^2 z^2 \TrS a c^2 x y^2 z^2 \TrS a^2
b c^2 x y^2 z^2 \TrS a c^3 x y^2 z^2 \TrS a c^4 x y^2 z^2 \TrS a^3
c y^3 z^2 \TrS a^3 c^2 y^3 z^2 \TrS a^2 b^2 c z^3 \TrS a b c^2 z^3
\TrS a b c^3 z^3 \TrS a b^2 c x z^3 \TrS b c^2 x z^3 \TrS a b^2
c^2 x z^3 \TrS b c^3 x z^3 \TrS b c^4 x z^3 \TrS a b^3 c x^2 z^3
\TrS b^2 c^2 x^2 z^3 \TrS b^2 c^3 x^2 z^3 \TrS a^2 b c y z^3 \TrS
a c^2 y z^3 \TrS a^2 b c^2 y z^3 \TrS a c^3 y z^3 \TrS a c^4 y z^3
\TrS a^2 b^2 c x y z^3 \TrS a b c^2 x y z^3 \TrS c^3 x y z^3 \TrS
a b c^3 x y z^3 \TrS c^4 x y z^3 \TrS c^5 x y z^3 \TrS a^3 b c y^2
z^3 \TrS a^2 c^2 y^2 z^3 \TrS a^2 c^3 y^2 z^3 \TrS a b c^2 z^4
\TrS a b c^3 z^4 \TrS a b^2 c^2 x z^4 \TrS b c^3 x z^4 \TrS b c^4
x z^4 \TrS a^2 b c^2 y z^4 \TrS a c^3 y z^4 \TrS a c^4 y z^4 \TrS
a b c^3 z^5 \TrS c^5 z^5$\pSkip

\noindent $\al_{22} = a^2 b^3 x^3 y^2 \TrS b c^2 x^3 y^2 \TrS a^2
b^3 x z^2
 $ \pSkip

\noindent $ \bt_{22} = a^3 b^2 x^2 y^3 \TrS a c^2 x^2 y^3 \TrS a^3
b^2 y z^2 $ \pSkip

By Equation  \Ref{eq:a1} there are two dominant monomials for
which: \pSkip \noindent $\al_{22} = a^2 b^3 x^3 y^2 \TrS b c^2 x^3
y^2 \TrS a^2 b^3 x z^2 \leq { \bf b c^2 x z^2 } $ \pSkip
\noindent $ \bt_{22} = a^3 b^2 x^2 y^3 \TrS a c^2 x^2 y^3 \TrS a^3
b^2 y z^2 \leq { \bf  a c^2 y z^2}$ \pSkip
Thus, $\al_{22}$ and $\bt_{22}$  are inessential for $f_{22}$ and
$g_{22}$ respectively, namely $f_{22} \eqR g_{22}$. \pSkip
\end{description}

Composing all together, $\pol{A} \eqR \pol{B}$ to complete the
proof of Lemma \ref{lem:zeroMA}.
\end{proof}

\begin{note*} To perform the symbolic computations in the proof of Lemma \ref{lem:zeroMA}
we were assisted by the mathematical software  Wolfram Mathematica
$6$.
\end{note*}
 \end{document}